\documentclass[preprint, 10pt]{elsarticle}

\usepackage{amssymb}
\usepackage{algorithm}
\usepackage{algorithmic}
\usepackage{tikz}
\usepackage{pgfplots}
\usepackage{microtype}
\usepackage{amsthm}
\usepackage{caption}

\usepackage{amsbsy}
\usepackage{amsmath}
\usepackage{bm}
\usepackage{color}
\usepackage{dcolumn}
\usepackage{tabularx}
\usepackage{bbm}

\newcommand{\vx}{{\mathbf x}}
\newcommand{\vy}{{\mathbf y}}
\newcommand{\vz}{{\mathbf z}}
\newcommand{\vv}{{\mathbf v}}
\newcommand{\vw}{{\mathbf w}}

\newcommand{\vp}{{\mathbf p}}
\newcommand{\vb}{{\mathbf b}}

\newif\ifnotesw \noteswtrue
\newcommand{\jeffrey}[1]{\ifnotesw  \textcolor[rgb]{0,0,1}{  $\spadesuit$Jeffrey:\ {\sf \bf \it #1}\ $\spadesuit$ }\fi}

\newcommand{\setalglineno}[1]{%
	\setcounter{ALC@line}{\numexpr#1-1}}
\newtheorem{theorem}{Theorem}[]

\newtheorem{lemma}[theorem]{Lemma}
\theoremstyle{definition}

\theoremstyle{remark}
\newtheorem*{remark}{Remark}

\journal{Journal of Computational and Applied Mathematics}

\begin{document}

\begin{frontmatter}



\title{Maximal Poisson-disk Sampling for Variable Resolution Conforming Delaunay Mesh Generation: Applications for Three-Dimensional Discrete Fracture Networks and the Surrounding Volume}


\author{Johannes Krotz\corref{cor1} \fnref{label1}}

\affiliation[label1]{organization={Department of Mathematics, Oregon State University},
            city={Corvallis},
            postcode={97330}, 
            state={Oregon},
            country={USA}}

\author{Matthew R. Sweeney\fnref{label2}}

\author{Carl W. Gable\fnref{label2}}

\author{Jeffrey D. Hyman\fnref{label2}}

\affiliation[label2]{organization={Computational Earth Science (EES-16), Earth and Environmental Sciences, Los Alamos National Laboratory},
            city={Los Alamos},
            postcode={87545}, 
            state={New Mexico},
            country={USA}}

\author{Juan M. Restrepo\fnref{label3}}

\affiliation[label3]{organization={Oak Ridge National Laboratory},
            city={Oak Ridge},
            postcode={37830}, 
            state={Tennessee},
            country={USA}}

\begin{abstract}
We propose a two-stage algorithm for generating Delaunay triangulations in 2D and Delaunay tetrahedra in 3D that employs near maximal Poisson-disk sampling.
The method generates a variable resolution mesh in $2$- and $3$-dimensions in linear run time.
The effectiveness of the algorithm is demonstrated by generating an unstructured 3D mesh on a discrete fracture network (DFN).  
Even though Poisson-disk sampling methods do not provide  triangulation quality bounds in more than two-dimensions, we found that low quality tetrahedra are infrequent enough and could be successfully removed to obtain high quality balanced  three-dimensional meshes with topologically acceptable tetrahedra.
\end{abstract}




\end{frontmatter}


\section{Introduction}

There are a number of methods used to model flow and the associated transport of chemical species in low-permeability fractured rock, such as shale and granite.
The most common are continuum models, which use effective medium parameters~\cite{gerke1993dual,lichtner2014modeling,neuman1988use,neuman2005trends,tsang1996tracer,zimmerman1993numerical} and discrete fracture network/matrix (DFN) models, where fractures and the networks they form are explicitly represented~\cite{cacas1990modeling,long1982porous,nordqvist1992variable}. 
In the DFN methodology, individual fractures are represented as planar $N-1$ dimensional objects embedded within an $N$ dimensional space.
Both conforming methods, where the mesh conform to intersections~\cite{mustapha2011efficient,mustapha2007new}, and non-conforming methods, which use more complex discretization schemes so the mesh does not need to be conforming ~\cite{berrone2013pde,erhel2009flow,pichot2010mixed,pichot2012generalized}, are currently in use.  
If the matrix surrounding the fracture network needs to be meshed, complications of mesh generation are compounded for both conforming and non-conforming methods~\cite{berre2018flow}.

While the explicit representation of fractures allows for DFN models to represent a wider range of transport phenomena and makes them a preferred choice, when linking network attributes to flow properties~\cite{hadgu2017comparative,hyman2018dispersion,hyman2019linking}, it also leads to unique and complex issues associated with mesh generation.

We propose a two stage algorithm that generates a conforming variable resolution triangular mesh on a three-dimensional discrete fracture network.  
The proposed algorithm uses maximal Poisson-disk  sampling to efficiently generate the mesh of each fracture with controlled mesh resolution. 
In a Poisson-disk sampling enforce a minimal distance between its nodes is enforced.
The first stage is based on the framework presented in \cite{fastvar}. It uses a rejection algorithm to generate an initial Poisson-disk sampling within a runtime that is linear in the number of vertices generated. The second phase is based on \cite{varrad} and adds further points to the sampling, maximizing its density without violating the restrictions of a Poisson-disk sampling.

Once a Poisson-disk sampling is generated, a conforming Delaunay algorithm~\cite{murphy2001point} is used to connect this point distribution where lines of intersection between fractures form a set of connected edges in the Delaunay triangulation of the network.
The time it takes to generate the samplings scales linearly with the number of nodes. While it is not guaranteed that the density of our Poisson-disk sampling is maximal, i.e. no further nodes can be added without violating the restrictions on distances between nodes, we show that in practice our samples are maximal enough to obtain high quality meshes.
We also present a three-dimensional version of the method that can be used to create a tetrahdron mesh of the volume surrounding the network that conforms to the fracture network. 

In section~\ref{sec:background}, we describe the challenges in the DFN mesh generation and the general properties of maximal Poisson-disk sampling. 
In section~\ref{sec:methods}, we provide a detailed explanation of our method, for both 2D fracture networks and 3D volume meshing. 
In section~\ref{sec:results}, we propose metrics to access the quality of the mesh and run times for both 2D and 3D demonstration examples.
In section~\ref{sec:conclusions}, we provide a few remarks.

\section{Background}\label{sec:background}

\subsection{Discrete Fracture Networks: Mesh Generation Background}

Due to the epistemic uncertainty associated with hydraulic and structural properties of subsurface fractured media, fracture network models are typically modeled probabilistically ~\cite{national2020characterization,national1996rock,neuman2005trends}.
In the DFN methodology, individual fractures are placed into the computational domain with locations, sizes, and orientations that are sampled from appropriate distributions based on field site characterizations.
The fractures form an interconnected network embedded within the porous medium. 
Each fracture  must be meshed for computation, so that the governing equations for flow and transport can be numerically integrated to simulate physical phenomena of interest.

Formally, each fracture in a DFN can be represented as a planar straight-line graph (PSLG) composed of a set of line segments that represent the boundary of the fracture and a set of line segments that represent where other fractures intersect it. 
Then each fracture can be described by a set of boundary points on the PSLG, denoted $\{p\}$, and a set of intersection lines $\{\ell_{i,j}\}$, where the subscripts $i$ and $j$ indicate that this line corresponds to the intersection between the $i$th and $j$th fractures. 
Once $\{p\}$ and $\{\ell_{i,j}\}$ are obtained for every fracture in the network, a point distribution covering each fracture must be generated.
If a conforming numerical scheme is used, then all cells of $\{\ell_{i,j}\}$ are discretized lines in the mesh which must coincide between intersecting fractures. 

So long as minimum feature size constrains are met, a conforming triangulation method, such as presented in ~\cite{murphy2001point}, can be implemented to connect the vertices such that all lines of intersection form a set of connected edges in a triangulation.

In general, one wants to properly resolve all relevant flow and transport properties of interest while minimizing the number of nodes in the mesh, and these two goals compete. 
While a uniform mesh resolution is fairly straightforward  and appropriate for Eulerian transport simulations, where spatially variable numerical diffusion can drastically affect fronts in the solute field~\cite{benson2017comparison}, the resulting mesh will have a large number of nodes;  computations will be more expensive. 
Variable mesh resolution can be appropriate for single-phase flow simulations or in particle tracking simulations where the spatially variable resolution does not adversely affect transport properties.
However, this variable mesh generation is more complex than its uniform counterpart. 
One of the principal complications of variable mesh resolution generation is creating a smooth transition of cell sizes. Absent that, jumps in the computed fields of interest and other numerical artifacts will occur. the starting point for the notion of mesh
quality would appear to be the analysis leading to
the minimum angle condition that the smallest angle
should be bounded away from zero. This perhaps
originated with Zlamal \cite{ZLAMAL} and is quoted by Strang
and Fix \cite{strangfix} together with a statement regarding how
poorly shaped triangles may have an effect on the condition
number of the linear algebra problem that must
be solved. This result was improved by Babuska and
Aziz \cite{azis}.
Most  methods for the generation of a conforming DFN mesh use a uniform point distribution on the networks and then modify the connectivity locally to conform to intersections~\cite{mustapha2007new,mustapha2011efficient}.
When using a conforming mesh, the numerical methods for resolving flow and transport in the network are typically simpler and have fewer degrees of freedom compared to non-conforming mesh methods~\cite{fumagalli2019conforming}.  
Similarly, almost all non-conforming numerical methods use a uniform resolution, but some  create variable resolutions across fractures (still uniform within a single plane) in an attempt to reduce the number of total nodes in the mesh~\cite{berrone2019parallel}.
A variable mesh resolution in non-conforming schemes could drastically reduce the number of nodes in the mesh while retaining the the ability to retain higher orders of accuracy. However it is rarely implemented due to the associated meshing complications~\cite{borio2021comparison}.

The generation of a variable resolution unstructured conforming mesh is quite rare, even with the advantages noted above. 
One technique in use is the Features Rejection Algorithm for Meshing ({\sc FRAM}) that addressed the issues associated with conforming DFN mesh creation by coupling it with network generation~\cite{hyman2014conforming}. 
Through this technique, {\sc FRAM} allows for the creation of a variable resolution mesh that smoothly coarsens away from intersections where pressure gradients are typically the highest in flow simulations. 
{\sc FRAM} has been implemented in the computational suite {\sc dfnWorks}~\cite{hyman2015dfnWorks}, which has been used to probe fundamental aspects of geophysical flows and transport in fractured media~\cite{hyman2020flow,hyman2019emergence,hyman2019matrix,kang2020anomalous,makedonska2016evaluating,sherman2020characterizing} as well as practical applications including hydraulic fracturing operations~\cite{hyman2017discontinuities,karra2015effect,lovell2018extracting}, inversion of micro-seismicity data for characterization of fracture properties~\cite{mudunuru2017sequential}, the long term storage of spent civilian nuclear fuel~\cite{hadgu2017comparative}, and geo-sequestration of carbon dioxide into depleted reservoirs~\cite{hyman2020characterizing}.

However, the implementation used is an iterative refinement method for point distribution,  which is very inefficient. 
To triangulate each polygon a `while` loop is was executed to apply a Rivara refinement algorithm to an initially coarse distribution based on the boundary set $\{p\}$.
If an edge in the mesh is greater than the current maximum edge length, a new point is added to the mesh at the midpoint of that edge to split it in two. 

In practice, the edge splitting is done using Rivara refinement~\cite{rivara1984algorithms,rivara1984mesh}.
The resulting field is then smoothed using Laplacian smoothing in combination with Lawson flipping~\cite{
khamayseh1996anisotropic}. 
This process is repeated until all edges met the assigned target edge length, which could be a spatially variable field based on the distance to  $\{\ell_{i,j}\}$, for example.
While the resulting mesh quality is quite good, the process is inefficient and cumbersome. 

The superior modeling qualities of variable resolutions can be made practical, if implementation complexities can be addressed. We do so here using 
a Poisson-disk 
sampling methodology where the final vertex distribution is directly created rather than iteratively derived. 
While the method was initially designed to specifically improve  {\sc FRAM}, we provide the details in a general format such that it can be implemented for any discrete fracture network methodology, including those that use both conforming and non-conforming flow and transport simulations. 
Details are given for Delaunay triangulations, which are of importance in many two-point flux finite volume solvers as they are used to generate the Voronoi control volumes on which these solvers compute.  
In the next section, we recount the properties of maximal Poisson-disk sampling that we used to design and implement this new method. 
Further we recount theoretical bounds on mesh gradation that ensure high-quality variable mesh resolutions.

\subsection{Maximal Poisson-disk Sampling}\label{sec:mpds}

Over the last few decades, there has been extensive research into mesh generation using maximal Poisson-disk sampling~\cite{MDS_summary}.
 Initially, it was driven by computer graphics aimed at removing artifacts in digitally-synthesized imagery~\cite{aliasing_3,aliasing_1,aliasing_4,aliasing_2}. 
Another application in computer graphics is the real time adjustment of the level of zoom in computer games~ \cite{adaptive_res}. 
When generating meshes from a sample, dense, yet cluster free samplings have provable high quality bounds \cite{MDS_good_trig2,MDS_good_trig3,MDS_good_trig1}. 
Similar quality bounds can be established for sphere-packings, whose radii are Lipschitz continuous with respect to their location \cite{Lipschitz1,Lipschitz2,Lipschitz3}.
Maximal, or almost maximal, Poisson-disk samplings fulfill all these conditions leading to high-quality meshes. 
It was shown in \cite{DFN_good1} and \cite{DFN_good2} that in fracture mechanics, where cracks propagate along mesh edges, meshes generated by a maximal Poisson-disk sampling generate more realistic cracks.
Traditionally, Poisson-disk sampling is generated with an expensive dart-throwing algorithm \cite{aliasing_3}.
These algorithms struggle to achieve maximality as the probability to select a free spot becomes decreasingly small. The algorithm in \cite{varrad} based on these dart-throwing algorithm is the first to guarantee maximality and reaches run times of $O(n\log(n))$ ($n:$ number of points sampled) by using a regular grid for acceleration and sampling from polygonal regions in its second phase to achieve maximality.
They report close to $O(n)$ performance in practice \cite{eff_mds,mds_simple,varrad}.
Prior to that  an algorithm not longer based on dart-throwing was proposed in \cite{Bridson2007FastPD}, which while not guaranteeing maximality, showed linear performance in the number of nodes sampled. 
Their algorithm was extended to variable radii \cite{fastvar}.

A maximal Poisson-disk sampling $X$ on a domain $\Omega\subseteq \mathbb{R}^d$ is a random selection of points $X=\{\vx_i\}_{i=1}^n$, that fulfills the following properties:
\begin{enumerate}
	\item \textit{empty disk property}: $$\forall i\neq j \in \{1,...,n\}: |\vx_i-\vx_j|>r.$$
	We will call $r$ the \textit{inhibition radius},
	\item \textit{maximality}:
	$$\Omega= \bigcup_{i=1}^n B_R(\vx_i),$$
	where $B_\varepsilon(\vx)=\{\vy\in \Omega: |\vx-\vy|<\vy\}$ is the open ball of radius $\varepsilon$ around $\vx$. $R$ will be called the \textit{coverage radius}. \cite{varrad}
	
\end{enumerate}
Intuitively, the \textit{empty disk property} says that  every center sample point of a $d$-dimensional ball or disk  does not contain any other points of the sampling. \textit{Maximality} implies that these balls cover the whole domain, i.e., there is no point $y\in \Omega$, that is not already contained in one of the balls around a point in the sample.

It is useful to generalize these definitions, such that both the \textit{inhibition} and the \textit{coverage radius} depend on the sampling points, i.e. $r=r(\vx_i,\vx_j)$ and $R=R(\vx_i,\vx_j)$ for all $\vx_i,\vx_j\in X$. We hereon refer to this construct as a variable radii  maximal Poisson-disk sampling, and we refer to a Poisson-disk sampling with constant radii as a fixed-radii maximal Poisson-disk sampling.\cite{varrad}

A common approach is to assign each point $\vx\in \Omega$ a positive radius $\rho(\vx)$ and have $r(\vx_i,\vx_j)$ be a function of $\rho(\vx_i)$ and $\rho(\vx_j)$. 
Natural choices for $r(\vx_i,\vx_j)$ are, for example, $\rho(\vx_i)$ or $\rho(\vx_j)$ for $i<j$, thereby determining the inhibition radius depending on the ordering on $X$. Order independent options include $\min(\rho(\vx_i),\rho(\vx_j)), \max(\rho(\vx_i),\rho(\vx_j))$ or $\rho(\vx_i)+\rho(\vx_j)$. The last of these options corresponds to a sphere packing \cite{varrad}. The coverage radius can, but does not have to be different from $\rho$.


The Delaunay triangulation of a sampling maximizes the smallest angle of its triangles among all triangulations of this sampling \cite{def_delaunay}. 
Since numerical errors in many applications tend to increase if these angles become smaller \cite{ZLAMAL}, Delaunay triangulations often are a triangulation of choice.
Moreover, the dual of the Delaunay triangulation is a Voronoi tessellation, which in a certain sense is optimal for two-point flux finite volume solvers~\cite{eymard2000finite}, that are commonly used in subsurface flow and transport simulators such as {\sc fehm}~\cite{zyvoloski2007fehm}, {\sc tough2}~\cite{pruess1999tough2}, and {\sc pflotran}~\cite{lichtner2015pflotran}. 
 In case of maximal Poisson-disk samplings we can go one step further and give a lower bound on these angles. 
 In what follows we estimate the bounds that apply to the sampling we generate on DFN in later sections.
 We provide a brief summary of the proofs found in \cite{varrad}, while highlighting the most important results we use. We then proceed with the new bounds.

\begin{lemma}
	\label{lemma:angle2d}
The smallest angle $\alpha$ in any triangle is grater than $\arcsin\left(\frac{r}{2R}\right)$, where $r$ is the length of the shortest edge and $R$ the radius of the circumcircle or \begin{align}
	\sin(\alpha)\ge \frac{r}{2R} \label{eq:smallest angle}
\end{align}
	\begin{proof}
		This is a direct corollary of the central angle theorem.
	\end{proof}
\end{lemma}
This Lemma allows us to give explicit bounds for maximal Poisson-disk samplings.
While we will focus entirely on inhibition radii given by $r(\vx_i,\vx_j)=\min(\rho(\vx_i),\rho(\vx_j))$, where $\rho(\vx)$ is some positive function, comparable results can be found for different $r(\vx_i,\vx_j)$ in a similar fashion. 
\begin{lemma}\label{lemma:2}
	Let $\varepsilon\ge 0$ and $\rho:\mathbb{R}^n\rightarrow \mathbb{R}$ ($n\ge 2$) be a positive Lipschitz continuous function with Lipschitz constant $L$ with $L\varepsilon < 1$.
	Let $X\subset$ be a variable maximal Poisson-disk sampling on the domain $\Omega\subset \mathbb{R}^n$ with inhibition radius $r(\vx,\vy)=\min(\rho(\vx),\rho(\vy))$  and coverage radius $R(\vx,\vy) \le (1+\varepsilon)r(\vx,\vy)$.($\varepsilon>0$)\\
	Let the triangle $\Delta$ be an arbitrary element of the Delaunay triangulation of $X$ ($n=2$) or an arbitrary 2-dimensional face of a cell of the Delaunay triangulation of $X$.  \\
	If the circumcenter of $\Delta$ is contained in $\Omega$, each angle $\alpha$ of $\Delta$ is greater or equal to $\arcsin\left(\frac{1-L-\varepsilon L}{2+2}\varepsilon\right)$ or 
	 \begin{align*}		 
	 \sin(\alpha)\ge \frac{1-L-\varepsilon L}{2+2\varepsilon}.
	\end{align*}
	\begin{proof}
		Let $\alpha$ be the smallest angle of $\Delta$ and $\vx,\vy\in X$ be the vertices of the shortest edge of $\Delta$, i.e. the vertices opposite to $\alpha$. 
		 Without loss of generality assume $\rho(\vx)\le \rho(\vy)$.
		 Since $X$ is a Poisson-disk sampling $|\vx-\vy|\ge \min\left(\rho(\vx),\rho(\vy)\right)=\rho(\vx)$.\\
		 Now let $\vz\in \Omega$ be the circumcenter of $\Delta$. Since $X$ is maximal, there exists $\vv\in X$ with $|\vz-\vv|\le R(\vz,\vv)\le (1+\varepsilon)\rho(\vz)$.  
Next we notice that, because $\Delta$ was retrieved from a Delaunay triangulation $\vv$ cannot be contained in the interior of $\Delta$'s circumcircle. 
Hence 	 \begin{align*}
	|\vz-\vx|\le |\vz-\vv|\le (1+\varepsilon)\rho(\vz)\le (1+\varepsilon)\left(\rho(\vx)+L|\vz-\vx| \right).
\end{align*}
Rearranging this inequality yields \begin{align*}
|\vz-\vx|\le \rho(\vx)\frac{1+\varepsilon}{1-L-\varepsilon L}.
\end{align*}
The result follows by applying Lemma \ref{lemma:angle2d} after noticing that $|\vx-y|$ is the length of the shortest edge and that $|z-\vx|$ is the radius of the circumcirle.
	\end{proof}
\end{lemma}
\begin{remark}
	Note that for $n>2$ the same result is true, if we assume the circumcenter of the $n$-simplex, whose face $\Delta$, is contained in $\Omega$ instead of the circumcenter of $\Delta$ itself. The proof is identical. 
\end{remark}
\begin{remark}
	While this result allows to control the quality of 2D-triangulations of maximal Poisson-disk samplings, it can also be used to gauge how close a given Poisson-disk sampling is to being maximal.
\end{remark}
The previous Lemma only gives us bounds on all triangles, if their circumcenters are contained in $\Omega$. 
The next two Lemmas will give sufficient conditions to guarantee exactly this as long as $\Omega$ is a polytope.
\begin{lemma}\label{lemma:circinside}
	Let $\Omega\subset\mathbb{R}^2$ be a polygonal region and $X$ a maximal Poisson-disk sampling containing all vertices of $\Omega$. Let the inhibition radius $r(\vx,\vy)$ be defined like in the previous lemma. 
	Further let the coverage radius of $X\cap \delta \Omega$ fulfill $R^\delta(\vx,\vy)<\frac{r(\vx,\vy)}{\sqrt{2}(1+L)}$, i.e. $|\vx-\vy|<{\frac{\sqrt{2}}{1+L}}r(\vx,\vy)$ for all $\vx,\vy\in \delta\Omega\cap X$, the circumcenter of all triangles in the Delaunay triangulation of $X$ are contained in $\bar{\Omega}$.
	\begin{proof}
		Suppose this claim is wrong. Then let $\Delta$ be a triangle in the Delaunay triangulation with circumcenter $z\notin \Omega$. For this to be possible the circumcircle needs to be cut in (at least) two pieces by $\delta\Omega$, separating $z$ and the vertices of $\Delta$. Since $\Delta$ is part of a Delaunay triangulation and all vertices of $\Omega$ are part of the sampling, this is done by (at least) one segment of a straight line, i.e. $\delta \Omega$ contains a secant of the circumcircle. \\
		Let $\vb_1,\vb_2\in \delta\Omega$ be the two boundary points closest to the circumcircle on either side of that line segment and let $B$ be the disk bounded by the circumcircle. Note that $\bar{B}\cap\Omega$ contains $\Delta$ and  is itself entirely contained in the disk of radius $\frac{1}{2}|\vb_1-\vb_2|<\frac{1}{2}{\frac{\sqrt{2}}{1+L}}r(\vb_1,\vb_2)$ around $\frac{1}{2}(\vb_1+\vb_2)$.\\
		Now let $\vx\notin \{\vb_1,\vb_2\}$ be a vertex of $\Delta$ and let $\vb\in\{\vb_1,\vb_2\}$ be the point of the two, that is closer to $\vx$. We already established that $\vx$ lies within the just mentioned ball around $\frac{1}{2}(\vb_1+\vb_2)$.
		Let $\vx_p$ be the projection of $\vx$ onto the line segment connecting $b_1$ and $b_2$. Then $|\vx-\vx_p|< \frac{1}{2}{\frac{\sqrt{2}}{1+L}}r(b_1,b_2)$, because $\vx$ lies within the circle of that radius, $|\vx_p-b|< \frac{1}{2}{\frac{\sqrt{2}}{1+L}}r(b_1,b_2)$, because $b$ is the closer of the two points ${b_1,b_2}$ and therefore  \begin{align}
		\label{eq:rhob}
		|\vx-\vb|=\sqrt{|\vx-\vx_p|^2+|\vb-\vx_p|^2}< \frac{r(\vb_1,\vb_2)}{1+L}\le \frac{\rho(\vb)}{1+L}\le \rho(\vb).\end{align}
		Since $|\vx-\vb|\ge \min(\rho(\vx),\rho(\vb))$ this implies $|\vx-\vb|\ge \rho(\vx)$.
		However assuming this and applying the Lipschitz condition on \eqref{eq:rhob} gives us 
		\begin{align*}
			|\vx-\vb| < \frac{\rho(\vb)}{1+L}\le \frac{1}{1+L}\left(\rho(\vx)+L|\vx-\vb|\right)\le \frac{1+L}{1+L}|\vx-\vb|,
		\end{align*}
		which is a contradiction.
		\end{proof}
\end{lemma}
\begin{remark}
	Lemma \ref{lemma:circinside} does generalize to higher dimensions. It is not very practical because
	 it is difficult to guarantee the bounds on $R^\delta$, if the boundary is more than 1-dimensional.
	However it is still possible to get some bounds on the radii of the circumcircles and then, using Lemma \ref{lemma:angle2d}, on the angles, if the distance of non-boundary nodes is greater than some lower bound $d>0$.\\

	In fact, using notation from the previous proof, let $\Delta$ again be an $n$-simplex with circumcenter outside of $\Omega$ and $\vx \notin \delta\Omega$ on of its nodes. Since the circumsphere of any simplex in a Delaunay triangulation does not contain any other nodes the radius of the intersection with $\delta \Omega$ is bounded by $R^\delta$.
	One can show with some simple geometry that this forces the radius of $\Delta$'s circumsphere $R$ to fulfill the following inequality \begin{align} \label{eq:outsidecenter}
		R^2\le (R-d)^2+\left(R^\delta\right)^2 \Rightarrow R\le \frac{d^2+\left(R^\delta\right)^2}{2d}\le \frac{\left(R^\delta\right)^2}{d}.
	\end{align}  
	If $R^\delta(\vx,\vy)<r(\vx,\vy)$ there is a lower bound on $d$, continuously depending on $R^\delta$, solely due to the fact, that we have a Poisson-disk sampling. If $R^\delta=R\delta(\vx_p,b)$ is any bigger, $d$ needs to be bounded artificially. This implies that the angle bounds change continuously, if the conditions for Lemma \ref{lemma:circinside} cannot be met they still can be relatively controlled by the choice of the artificial bound on $d$.
\end{remark}

Under the conditions of the previous Lemmas the simplices of the Delaunay triangulation is guaranteed to only have well-behaved triangular faces. In three or more dimensions however this does not imply that the simplices themselves are well-behaved.
It is still possible for a Delaunay triangulation to contain slivers for example, that is tetrahedra whose 4 nodes are all positioned approximately on the equator of their circumsphere. In \cite{exudation} slivers are characterized as tetrahedra, whose nodes are all close to a plane and whose orthogonal projection onto that plane is a quadrilateral. In \cite{bern} slivers are equivalently classified as tetrahedra with a dihedral angle close to $180^\circ$ containing their own circumcenter. 
Slivers can have all their faces equilateral triangles, yet have  dihedral angles that are arbitrarily small, causing numerical errors to blow up. \\

While slivers cannot be entirely avoided, one can show that if the nodes $\vx,\vy,\vz,\vw$ of a maximal Poisson-disk sampling form a sliver, the distance between $\vw$ and the plane spanned by $\vx,\vy,\vz$ needs to be very small \cite{exudation}.
This allows us to avoid slivers within certain planes, by first generating a 2D sampling in these planes and then enforcing a minimal distance between the plane and further nodes in the 3D sampling.  We use this to avoid slivers around the DFN and the faces of the surrounding matrix.
This also causes slivers to be rather scarce in a 3D maximal Poisson-disk sampling as given any three nodes the vast majority of possible positions of a fourth node do not produce a sliver.   
This scarcity of slivers in a sampling makes it quite likely that if nodes of slivers are removed and resampled the resulting triangulation will have less slivers than the previous one. This opens the door of a rejection-style algorithm to be successful in improving the overall quality of a triangulation.\\


\section{Methods}\label{sec:methods}
Our proposed method for mesh generation is broken into three primary steps. 
First, we generate a 2D point distribution upon each fracture in the DFN. 
After merging these samples and removing conflicts with regards to the empty disk property, we generate a 3D-Poisson disk sampling on the surrounding matrix by adding points wherever maximality allows it. 
Finally, in an attempt to remove slivers, we remove their nodes and randomly replace them until no slivers remain.

\subsection{2D Sampling Method}
We generate 2D Poisson-disk samplings in a successive manner using a rejection method. 
This method can be preformed on every fracture in the network independent of the other fractures. 
(Details can be found in~\cite{hyman2014conforming}.)
In each step a new candidate is generated, and if it does not break the empty disk-property with any of the already accepted nodes, it is accepted. 
For the sampling in two dimensions, we use a variable inhibition radius that increases linearly based on the distance to the closest intersection of the DFN. 

In particular, we reject a candidate node $\vy$, if there is an already excepted node $\vx$ such that the condition \begin{align}
|\vy-\vx|\ge r(\vx,\vy) = \min(\rho(\vx),\rho(\vy)) \label{eq:empty_disk}
\end{align}
is violated.
In this equation $\rho(\vx)$ as a piecewise linear function given by
\begin{align}
\rho(\vx)=\rho(D(\vx))=\left\{\begin{matrix}
\frac{H}{2} &\text{for} & D(\vx)\le FH\\
A(D(\vx)-FH)+\frac{H}{2} &\text{for}& FH\le D(\vx) \le (R+F)H \\
\text{else} \ \ (AR+\frac{1}{2})H 
\end{matrix}\right. .\label{eq:rho2d}
\end{align}
Here $D(\vx)$ is the Euclidean distance between $\vx$ and the closest intersection. 
$H,A,R$ and $F$ are parameters, that determine the global minimal distance between two nodes ($H/2$), the range around an intersection on which the local inhibition radius remains at its minimum ($FH$), the global maximal inhibition radius ($ARH+H/2$), and the slope at which the inhibition radius grow with $D(\vx)$ ($A$).
Since $\rho(D)$ is piecewise linear, it is a Lipschitz-function with Lipschitz-constant $A$. 

\begin{figure}
	\centering
	\includegraphics[height=0.4\linewidth]{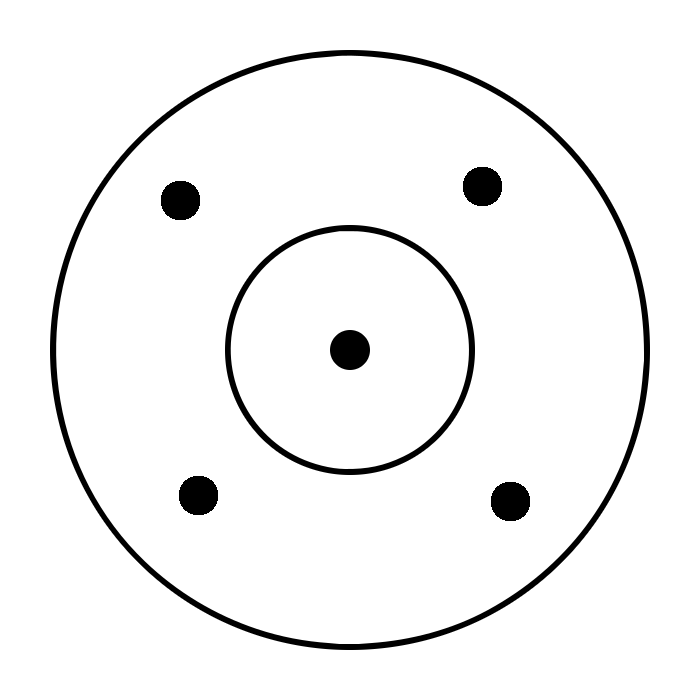}
	\caption{Visualisation of a single sampling step. Current node at center, new candidates in annulus (k=4). Inner circle is bounded by the inhibition radius of the current node. Outer circle is bounded by maximal distance a node could be away from the current if the Poisson-disk sampling was maximal. }
	\label{fig:annulus}
\end{figure}

If the sampling has a coverage radius $R(\vx,\vy)\le (1+\varepsilon)r(\vx,\vy)$ for some $\varepsilon>0$ the conditions of \eqref{lemma:2} hold. 
To satisfy the conditions of \eqref{lemma:circinside} as well and thereby ensure angle bounds on all triangles in a Delaunay triangulation we first sample along the boundary, enforcing a maximal distance of $\frac{r(\vx,\vy)}{\sqrt{2}(1+L)}$ between boundary nodes.
As shown in \cite{Bridson2007FastPD} and \cite{fastvar}, we generate new candidates for our sampling randomly on an annulus  around an already accepted node. 
This is illustrated in Figure \ref{fig:annulus}.  
The inner radius of this annulus is determined by the minimal distance another node could have to the center node, while still preserving the empty disk property, whereas the outer radius is determined by the maximal distance a node could have to the center in a maximal sampling.  
For our choice of inhibition radius, assuming the same radius as coverage radius, these distances can be made out to be $r_{in}(\vx)=\frac{\rho(\vx)}{1+A}$  and $r_{out}(\vx)=\frac{2\rho(\vx)}{1-A}$.

%
%
%

\begin{figure}
	\centering
(a)	\includegraphics[height=0.4\linewidth]{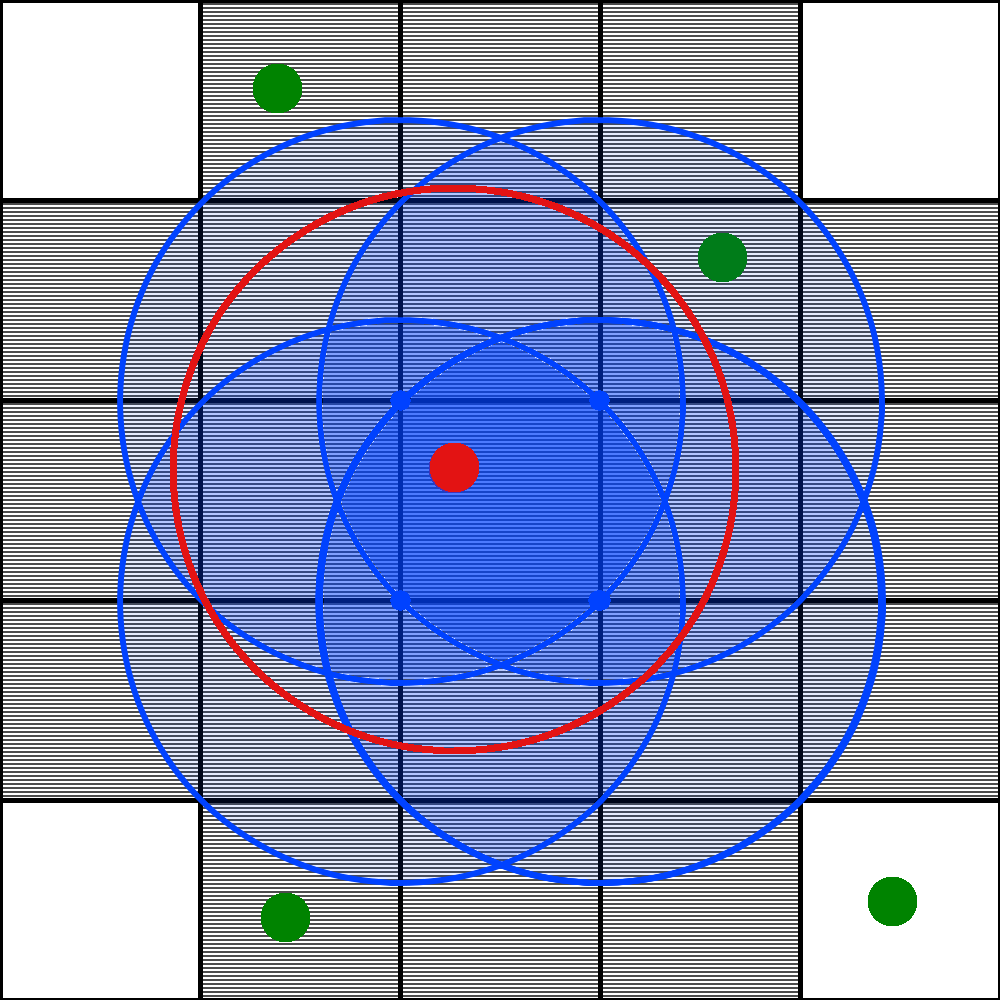} \hspace{0.5cm}
(b)		\includegraphics[height=0.4\linewidth]{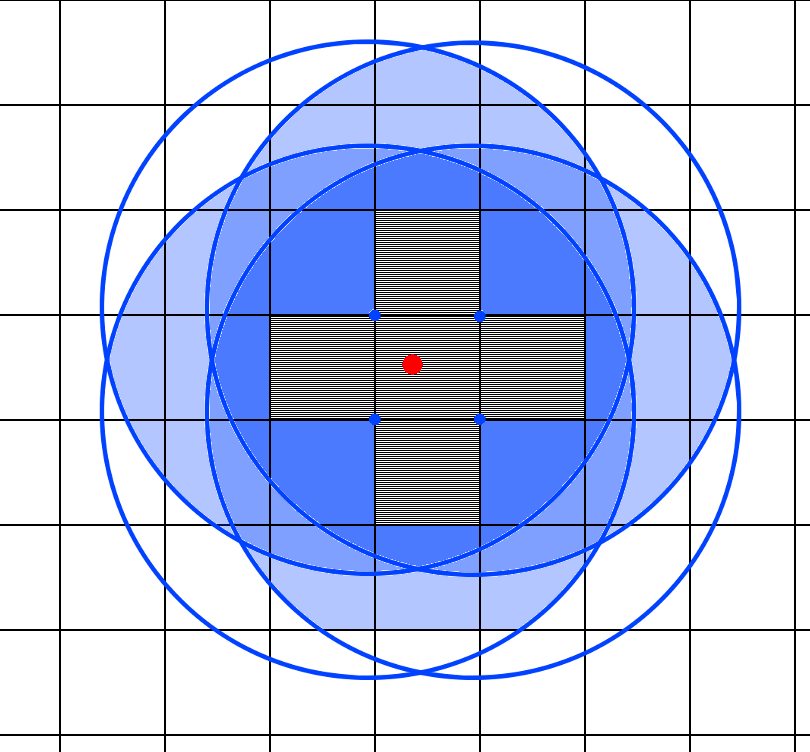}
	\caption{Visualisation of how the grid is used to find possibly conflicting nodes. New candidate in red, already accepted nodes in green, cells that can contain conflicting nodes in grey. Red circle shows the inhibition radius of the candidate, blue circles show furthest cells a node in the center cell could conflict with.}
	\label{fig:grid}
\end{figure}

\begin{figure}
{\newpage
\hspace{-2cm}
\begin{minipage}[c]{0.45\textwidth}
	\vspace{-2cm}
	\begin{tikzpicture}
		\node (initial) at (.4,0) {\includegraphics[width=.2\textheight]{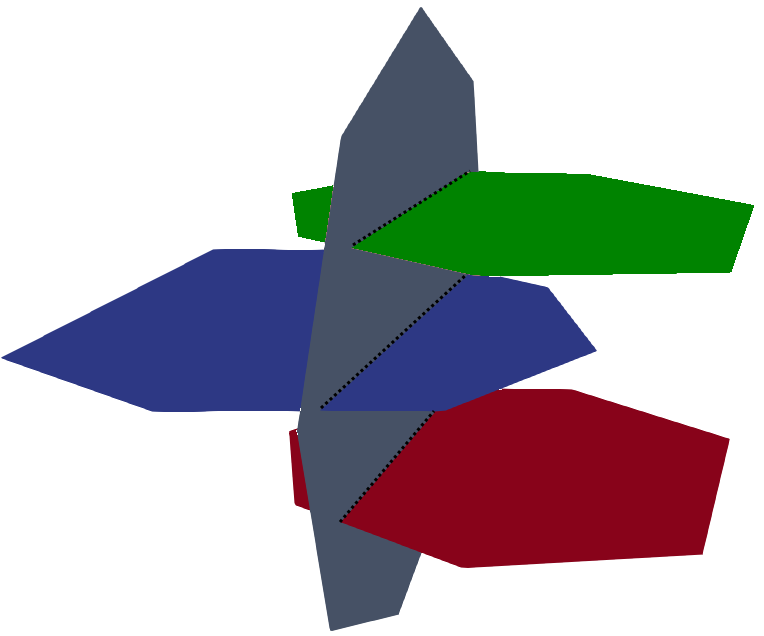}};
		\node (decomp) at (0, -4.5) {\includegraphics[width=.3\textheight]{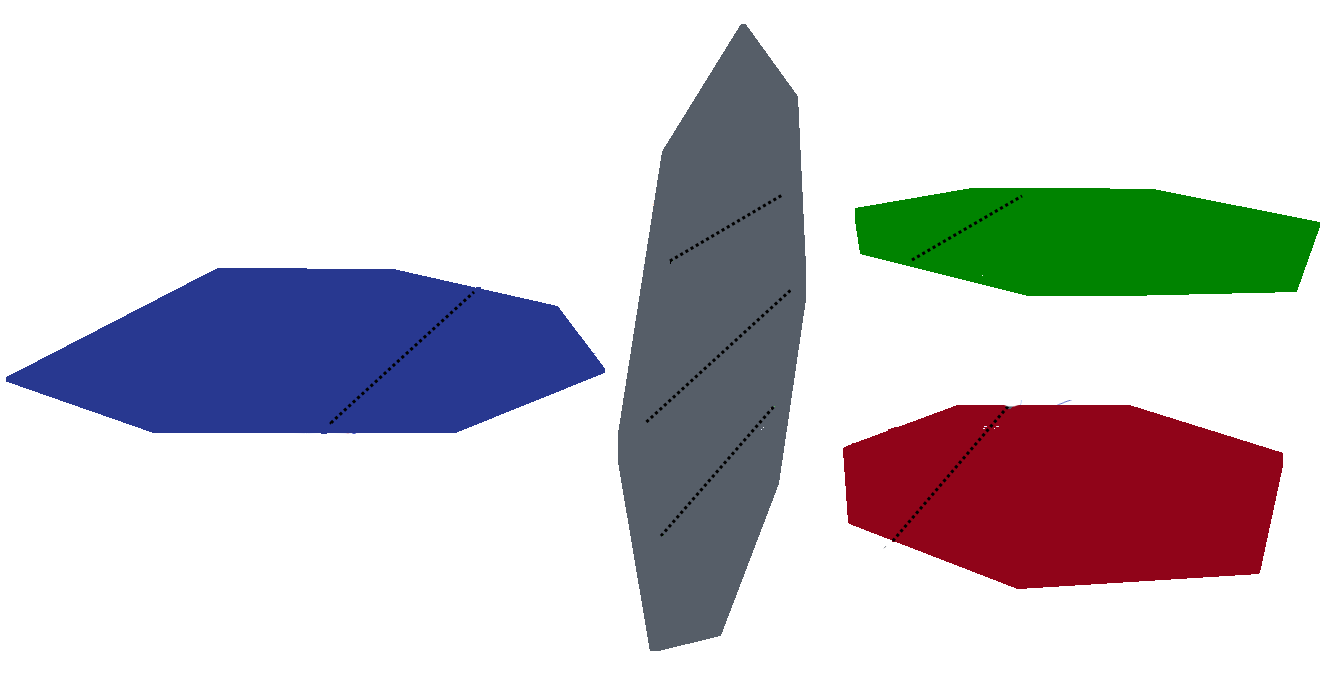}};
		\node (decomp) at (-3, -9) {\includegraphics[width=.3\textwidth]{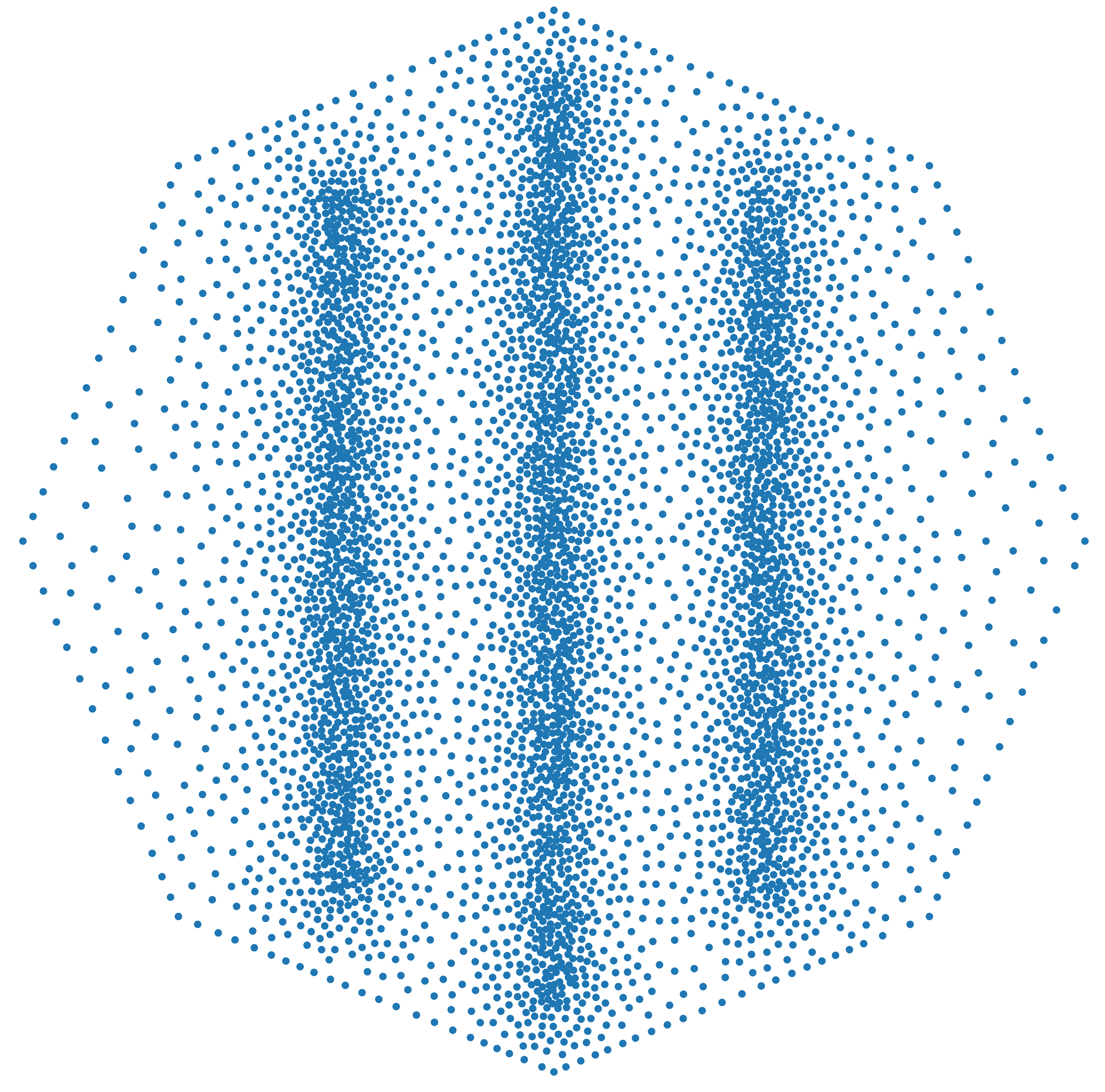}};
		\node (decomp) at (3, -9) {\includegraphics[width=.3\textwidth]{graphics/individual}};
		\node (decomp) at (0, -9) {\includegraphics[width=.3\textwidth]{graphics/individual}};
		\node (decomp) at (0.4, -13) {\includegraphics[width=.6\textwidth]{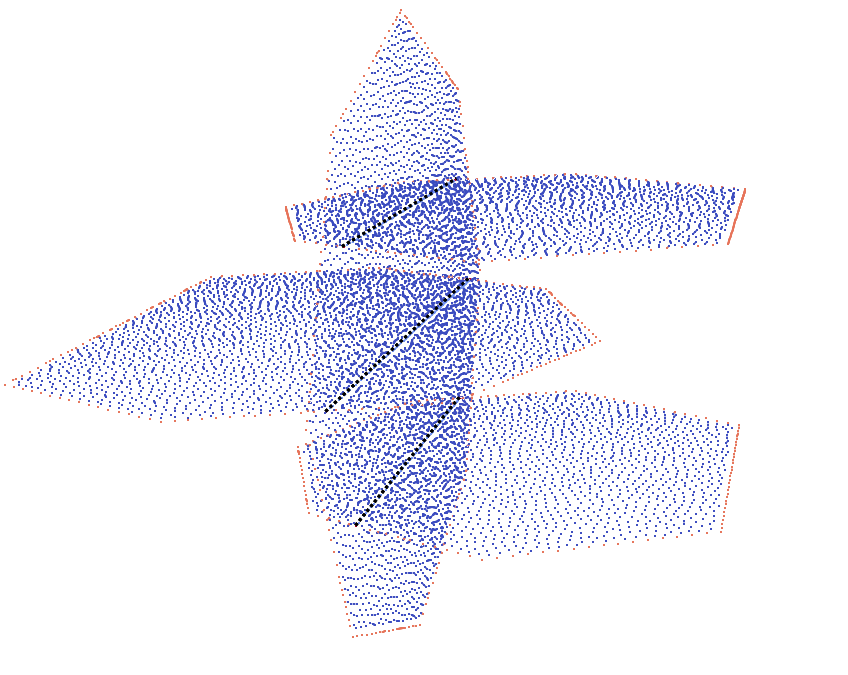}};
		\node (decomp) at (-1.2, -18) {\includegraphics[width=.7\textwidth]{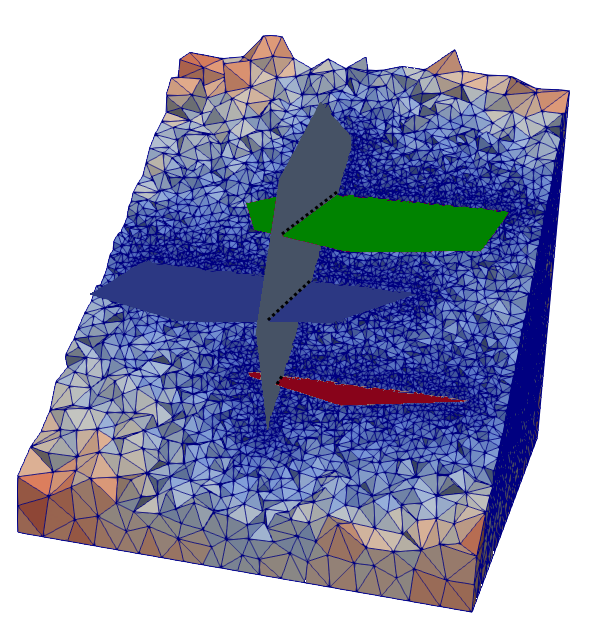}};

		\draw[->,black,ultra thick] (-1.,-2) -- (-1.5,-3)
		node[midway,fill=white]{\hspace{.5cm }};
		\draw[->,black,ultra thick] (1.4,-2) -- (1.9,-3)
		node[midway,fill=white]{      	\hspace{.5cm }};
		\node[fill=white] () at (.4,-2.5) { 1. decompose into individual fractures};

		\draw[->,black,ultra thick] (-2.2,-6.5) -- (-2.7,-7.5)
		node[midway,fill=white]{\hspace{.5cm }\vspace{1cm}};
		\draw[->,black,ultra thick] (-.0,-6.5) -- (-0,-7.5)
		node[midway,fill=white]{\hspace{.5cm }\vspace{1cm}};
		\draw[->,black,ultra thick] (2.,-6.5) -- (2.5,-7.5)
		node[midway,fill=white]{\hspace{.5cm }\vspace{1cm}};
		\node[fill=white] () at (0,-7) {2.sample individually (possibly parallel)};

		\draw[<-,black,ultra thick] (-1.,-11.5) -- (-1.5,-10.5)
		node[midway,fill=white]{\hspace{.5cm }};
		\draw[<-,black,ultra thick] (1.4,-11.5) -- (1.9,-10.5)
		node[midway,fill=white]{      	\hspace{.5cm }};
		\node[fill=white] () at (.4,-11.) {3.reassemble and eliminate conflicts};

		\draw[->,black,ultra thick] (0.,-15) -- (-0.3,-16)
		node[midway,fill=white]{\hspace{.5cm }};
		\node[fill=white] () at (-.2,-15.5) {4.sample matrix};

		\draw[-,black, ultra thick] (1.6,-19.5) arc(-90:90:2) -- (1.6,-15.5)  ;
		\draw[->,black, ultra thick] (1.6,-15.5)--(1.55,-15.5);
		\node[text width=2.5cm,align=center] at (2,-17.5) {5. remove slivers  and repeat};
		\end{tikzpicture}

\end{minipage} \hspace{3cm}
\begin{minipage}[c]{0.45\textwidth}
	\vspace{-4cm}
	
	\begin{tikzpicture}
	\node (init) at (5,0) {\includegraphics[width=0.7\linewidth]{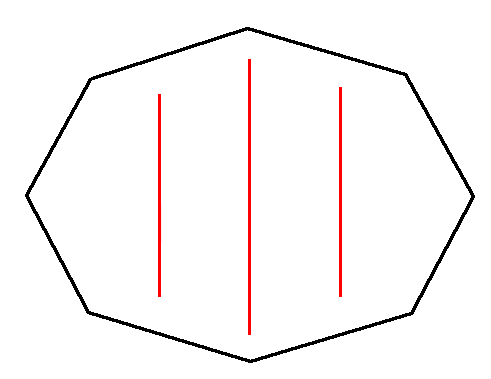}};
	\node (init) at (3,-4) {\includegraphics[width=0.7\linewidth]{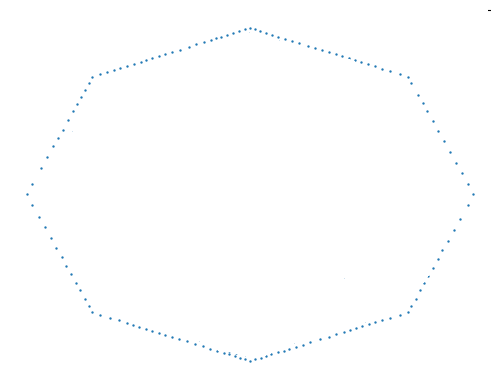}};
	\node (init) at (5,-8) {\includegraphics[width=0.7\linewidth]{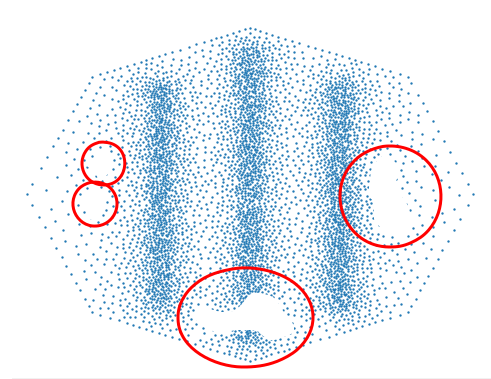}};
	\node (init) at (3,-12) {\includegraphics[width=0.7\linewidth]{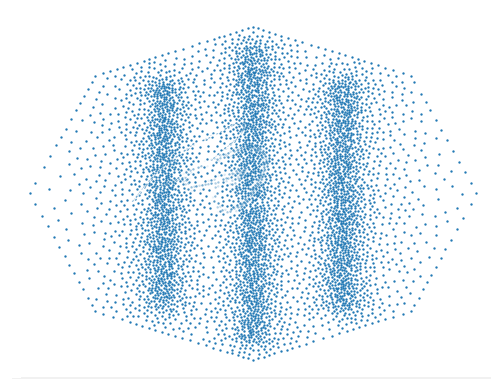}};
	
	\draw[->,black,ultra thick] (4,-1.5) -- (3.7,-2.5);
	\draw[->,black,ultra thick] (3.7,-5.5) -- (4,-6.5);
	\draw[->,black,ultra thick] (4,-9.5) -- (3.7,-10.5);
	
	\node[text width=5cm] () at (4,2) {1.Identify boundary and intersections for faster acces};
	\node[fill=white] () at (3.85 ,-2.) {2.Sample along boundary};
	\node[fill=white] () at (3.85 ,-6) {3.Use as seed for sampling};
	\node[fill=white, text width = 3cm] () at (6.85 ,-11.5) {4.When algorithm halts find undersampled regions and sample nodes in there. Repeat from step 3.};	
	\end{tikzpicture}
	\vspace{0.5cm}

\caption{\label{fig:workflow} Overview of workflow between creation of DFN and final mesh (left) and overview of workflow during 2D-sampling (right). }	
\end{minipage}}
\end{figure}

We will now go over the individual steps of the 2D algorithm. These steps can also be found in the pseudocode Algorithm 1 in Section \ref{sec:pseudo} and are illustrated in figure \ref{fig:workflow}. The necessary notation to read the pseudocode is found in the table at the start of the same section. 
In line \ref{line:seed} of that code a 1D Poisson-disk sampling along the boundary of the polygon is generated as a seed to start the algorithm. 
We continue to sample $k$ new candidate nodes at a time (line \ref{line:k_samp}) around each already accepted node and determine whether they get accepted or not (lines \ref{line:rej} through \ref{line:reject_end}). 
$k$ is a positive integer and a user-defined parameter of the algorithm. 
If all $k$ candidates around a node are rejected, we move on to the next already accepted node (line \ref{line:next_node}). The algorithm terminates for the first time as soon as every accepted node was the sampling center once(line \ref{line:end}).
Following \cite{Bridson2007FastPD} and \cite{fastvar}, we use cell-lists to find nodes around a candidate that could potentially cause this candidate to violate the empty-disk property, as depicted in Figure \ref{fig:grid}(a). The size of these cells is chosen to contain at most one node. This allows us to disregard distance calculation with nodes beyond a certain cutoff and therefore allows us to achieve linear run times in the number of generated nodes(line \ref{line:yinNp}). 
However, unlike the previously mentioned algorithm we do not only label cells containing particles as occupied, but also cells that are too close to an accepted node to contain a particle.  In particular, if a candidate $\vx$ lies in a cell $C$ and any other cell $D$ with $diam(C\cup D)\le r_{in}(\vx)$ is occupied, $\vx$ can be rejected right away as it conflicts with the node in $D$ (line \ref{line:rej}). On the other hand, if $dist(C,D)>\rho(\vx)$, there is no need to calculate the distance between $\vx$ and any potential element of $D$, as they can never violate the empty disk-property. An example of that is shown in Figure \ref{fig:grid}(b).
We use this to our advantage in two ways:
First, it allows us to reject many candidates without calculating any distances to nearby nodes, which particularly for large values of $k$ gives a respectable speed up compared to the original algorithm; 
second, unmarked cells are easy to find  and contain at least  some space for another node, allowing us to find undersampled regions after the algorithm terminated(line \ref{line:emptycell}). 
We fill these holes in the sample by generating random candidates within these unmarked cells(line \ref{line:rand_new}). 
The main algorithm is then restarted from these newly added nodes till it terminates again(line \ref{line:rerun}). While this process can be repeated several times just a single resampling already increases the quality of the sampling tremendously. 

Once the point distribution is created, the conforming Delaunay triangulation method of~\cite{murphy2001point} is used to create the final mesh on the fracture. 
In order for a conforming Delaunay triangulation which preserves the lines of fracture intersections as a set of triangle edges is created, it is sufficient that the circumscribed circle of each segment of the discretized line of intersection be empty of any other node in the point distribution prior to connecting the mesh.
To achieve this condition, any node within the circumscribed circle of each segment of the discretized lines of intersection is removed from the point distribution. 
Next, a two-dimensional unconstrained Delaunay triangulation algorithm is used to connect this node set. 
Because of the construction method, i.e., empty regions around the lines of intersection, the line segments that represent lines of fracture intersection must emerge in the triangulation and the Delaunay triangulation will conform to all of the fracture intersection line segments. 
Once every fracture polygon is triangulated, they are all joined together into a unified triangulated fracture network. 

\subsection{3D sampling method}
The sampling in 3D works very similar to its 2D counterpart. However, new candidates are generated on a spherical shell around accepted nodes instead of on an annulus. The 3D variant of $\rho(\vx)$  given by
\begin{align}
\rho(\vx) = \rho(D(\vx))=\left\{\begin{matrix}
\rho_{2}(\vx_p)  &\text{ for }&  D(\vx)\le F\rho_2(\vx_p)\\\\
A(D(\vx)-Fr_{2D}(\vx_p))+\frac{H}{2} & \text{ for }& F\rho_{2D}(\vx_p)\le D(\vx) \le \frac{\rho_{max}-r_{2}(\vx_p)}{A} \\\\
\text{ else } \ \ \rho_{max}& & 
\end{matrix}\right. .\label{eq:rho3d}
\end{align}
$\vx_p$ the fracture point closest to $\vx$  and $\rho_2(\vx_p)$ is its 2D inhibition radius on the fracture. $D(\vx)$ is the distance between $\vx$ and $\vx_p$.
Like in 2D, this is a piecewise linear function in $D(\vx)$, which is constant, if within a  distance of $\rho_2(\vx_p)F$ ($F$ a parameter) and then increases linearly with a slope of $A$ until the maximal inhibition radius of $\rho_{max}$ is reached. In addition to rejecting all candidates $y$ for which \eqref{eq:empty_disk} is violated, we also reject a candidate $\vx$, if it is within a distance of $\rho(\vx)/2$ to a boundary or fracture. This both prevents slivers from having three nodes located on a single fracture or the boundary of the matrix and limits the circumradius of tetrahedra with circumcenter outside of the matrix (lemma \ref{lemma:circinside} and subsequent remark).

A pseudocode of how the 3D-sampling is run from here can be found in Section \ref{sec:pseudo} in Algorithm 2 . 
The necessary notation is listed in the table at the start of that section.
As the first sampling process is essentially identical to the 2D version, we will explain the differences in the initialization and the resampling.
At the start the nodes are initialized through a Poisson-disk sampling on the boundary of the 3D matrix and the sampling on the DFN generated by the 2D algorithm(line \ref{line:seed3d}).
Neighbor cells can still be used in the same way as in 2D to speed up the rejection of candidates. 
Unlike in 2D, a maximal Poisson-disk sampling does not guarantee sliver-free triangulation, which is why we do not use the cell lists to find undersampled cells in 3D. 
Instead, once the algorithm terminates, the resulting sampling is triangulated (line:\ref{line:triangulate}), slivers identified (line:\ref{line:sliver}), and 2 nodes of every sliver (with a preference for nodes, that are neither on a boundary or a fracture) removed (line \ref{line:rem_nodes}). While the definition of a sliver given earlier in section \ref{sec:mpds} allows for a bit of leeway in what is considered a small or large dihedral angle, in practice we successfully replaced tetrahedra with dihedral angles outside of $[8^\circ,170^\circ]$ and aspect ratios bigger than $0.2$.
Then the algorithm is restarted with the remaining nodes as seed(line \ref{line:rerun3d}). 
This process is repeated till a sliver-free sampling is obtained(line \ref{line:term3d}).
With this approach we have been able to obtain triangulations with no elements of dihedral angles of less than  $8^\circ$ (presented in next sections).
The method for forming the conforming mesh is similar to that for the two-dimensional case, but spheres around triangle cells of the fracture planes are excavated. 
Additional details are found in~\cite{hyman2021flow}.

%

\subsection{Workflow overview}
The workflow is depicted in Figure \ref{fig:workflow} and contains the following high-level steps: 
(1) generation of a DFN using dfnWorks \cite{dfn}, 
(2) decomposition of DFN into individual polygons, 
(3) generation of 2D-variable-radii Poisson-disk samplings on each individual polygon using algorithm 1\ref{alg:2d} below, 
(4) construct a conforming Deluanay triangulation as previously described, 
(5) merge individual fracture meshes into a sampling on the original DFN, removing conflicting nodes along intersections, 
(6) generating a conforming 3D variable radii Poisson-disk sampling of the surrounding matrix of the DFN using the 2D samplings as seed according to Algorithm \ref{alg:3d}, 
(7) triangulate sampling, identify low-quality tetrahedra and remove 2 of their nodes that are not located on the original DFN, 
(8) repeat steps 5 and 6 with the remaining nodes as seed until no more low-quality tetrahedra remain \cite{monte-carlo}. 
Replacing step (8) with more traditional ways of sliver-removal like perturbation \cite{Perturbation} or exudation \cite{exudation} can break the empty disk property of the sampling. 
\subsection{Pseudocode for the 2D and 3D sampling algorithms\label{sec:pseudo}}
\begin{tabularx}{\textwidth}{l X}
		\multicolumn{2}{l}{ \textbf{Notation for Pseudocodes:}}\\ 	\multicolumn{2}{l}{ \hrulefill}\\
	\multicolumn{2}{l}{ \textbf{Input:}}   \\
	\textbullet $D^3\subset\mathbb{R}^3$: & cubical domain ($\ast$) \\   
	\textbullet $DFN\subset D^3$: & generated by DFNWorks ($\ast$) \\
	\textbullet	$F_l \subset \mathbb{R}^3$  :& $l$-th fracture of the DFN  \\  \textbullet  $q_{l,m}^{(1)}$ and $q_{l,m}^{(1)}$:& endpoints of intersection between fractures $F_l$ and $F_m$   \\  
	\multicolumn{2}{l}{\textbf{User defined parameters: }}   \\
\textbullet  $H/2$: &minimal distance between nodes \\
 \textbullet  $F$ : & $HF$ is range of constant density around intersections \\
 \textbullet  $R$:& $ARH+H/2$ is maximal distance between nodes  \\
 \textbullet $A$:& max. slope of inhibition radius  \\
  \textbullet$k$:& number of concurrently sampled candidates\\
 \multicolumn{2}{l}{\textbf{Additional notation:}}   \\
 \textbullet $G$:& square cells covering $F_l$ with  $diam(g)\le H/2$ for all $g\in G$.\\
 \textbullet $\rho(\vx)$: & as defined in equation \eqref{eq:rho2d}(2D) or \eqref{eq:rho3d} (3D) \\
 \textbullet $r(\vx,\vy)$: &  inhibition radius $\min(\rho(\vx),\rho(\vy))$ \\
 \textbullet $R(\vx,\vy)$: & coverage radius \\
\textbullet$C(\vx) \in G$:& grid cell containing the point $x$.\\
 \textbullet $N^+(\vx)$ :& $\{g\in G: dist(C(\vx),g)\le \rho(\vx)\}$: cells that can contain points $y$ with $|\vx-y|\le r(\vx,y)$\\
 \textbullet $N^-(\vx)$: &$\{g\in G: diam(g\cup C(\vx))\le \frac{\rho(\vx)}{1+A} \}$: cells, where for all their points $y$ $|\vx-y|\le r(\vx,y)$\\
 \textbullet $G_{occ}$& $\bigcup\limits_{\vx\in X} N^-(\vx)$: cells on which $X$ is already maximal.\\
 \textbullet $\mathcal{T}(X)$: & Delaunay triangulation of $X$ ($\ast$)\\
 \textbf{Output:} &  \\
 \textbullet  $X$ :& Poisson-disk sampling on the $l$-th fracture \\ 
 	\multicolumn{2}{l}{\hrulefill}\\
 	$(\ast)$: 3D only 
\end{tabularx}

\begin{algorithm}[H]	
	\caption{1}{\textbf{2D Poisson-disk sampling}\label{alg:2d}}

	\begin{algorithmic}[1]		
		\STATE\textbf{Initializing:}
		\STATE $X \subset F_l$ \hfill $\triangleright$ Generate a 1D Poisson-disk sampling with $R(\vx,y)\le \frac{r(\vx,y)}{\sqrt{2}(1+L)}$\\ \hfill along boundary $\delta F_l$ as seed. \label{line:seed}
		\FOR{$\vx\in X$}
		\STATE $G_{occ}\leftarrow G_{occ}\cup N^-(\vx)$ \hfill $\triangleright$ initialize occupied cells
		\ENDFOR
		 
		 \hspace{0cm}
		\STATE\textbf{Sampling:}
		\STATE $i\leftarrow 1$ \hfill $\triangleright$ Start sampling at first accepted node.
		\STATE $N\leftarrow |X|$ 	\hfill $\triangleright$ will increase as more nodes are accepted
		\WHILE{$i\le N$}  \label{line:sampling}
		
		\REPEAT
		\FOR{$j \in \{1,...,k\}$} \label{line:gen_candidates}
		\STATE	$\vp_j \in F_l$ \hfill $\triangleright$ generate $k$  new candidate nodes on the annulus around $\vx_i$\label{line:k_samp}
		\IF{$C(\vp_j) \in G_{occ}$} \label{line:rej}
		\STATE reject $\vp_j$ 	\hfill $\triangleright$ Cell already blocked by existing node's inhibtion radius
		\ELSE 
		\FOR{$\vy\in N^+(\vp_j)$}\label{line:yinNp}
		\IF{$|\vp_j-\vy|<r(\vp_j,\vy)$}
		\STATE reject $\vp_j$ 	\hfill $\triangleright$ empty disk property violated
		\ENDIF
		\ENDFOR
		\ENDIF \label{line:reject_end}
		\IF{$p_j$ was not rejected} \label{line:accept_start}
		\STATE  $X \leftarrow X\cup \{p_j\}$\hfill $\triangleright$ accept $\vp_j$ and add it to the sampling
		\STATE $G_{occ}\leftarrow G_{occ} \cup N^-(\vp_j)$ \hfill $\triangleright$ update occupied cells	
		\STATE $N\leftarrow N+1$ 	\hfill $\triangleright$ ensures sampling around newly accepted \label{line:accept_end} nodes				
		\ENDIF		
		\ENDFOR 
		\UNTIL{All $k$ of the $\vp_j$ are rejected}\label{line:all_k_rejected}
		\STATE $i\leftarrow i+1$ \hfill $\triangleright$  start sampling around next accepted node \label{line:next_node}
		\ENDWHILE \hfill $\triangleright$ terminate here or start resampling \label{line:end}
	\end{algorithmic}
\end{algorithm}

\begin{algorithm}[H]
\captionsetup{labelformat=empty}
\caption{}{\textbf{Continuation of Algorithm 1 (2D Resampling)}}

	\begin{algorithmic}[1]	
		\setalglineno{31}	
		\STATE  \textbf{Resampling:}(optional: algorithm terminates, if no resampling is required.)
		
		\FOR{$C\in G \setminus G_{occ}$}\label{line:emptycell}
		\STATE $\vp \in C$  \hfill $\triangleright$ generate a random candidate on each cell \label{line:rand_new}
		\FOR{$\vy\in N^+(\vp)$}
		\IF{$|\vp-\vy|<r(\vp,\vy)$}	
		\STATE reject $\vp$ \hfill $\triangleright$ empty disk property violated
		\ENDIF
		\ENDFOR		
		\IF{$\vp$ was not rejected} 
		\STATE  $X \leftarrow X\cup \{\vp\}$\hfill $\triangleright$ accept $\vp$ and add it to the sampling
		\STATE $G_{occ}\leftarrow G_{occ} \cup N^-(\vp)$ \hfill $\triangleright$ update occupied cells 	
		\STATE $N\leftarrow N+1$

		\ENDIF

		\ENDFOR
		\STATE
		\STATE \textbf{Rerun algorithm again from line \ref{line:sampling}} ($i$ is not reset.)\label{line:rerun}
	\end{algorithmic}

\end{algorithm}

\begin{algorithm}[H]	
	\captionsetup[FLOAT_TYPE]{labelformat=simple}
	\caption{2}{\textbf{3D Poisson-disk sampling + Resampling}\label{alg:3d}}

	\begin{algorithmic}[1]
		\STATE \textbf{Initializing:}
		\STATE $X\subset D^3$ \hfill $\triangleright$ Use Algorithm 1 to generate a Poisson-disk sampling on $\delta D^3$ and the DFN by using Algorithm 1 (remove conflicting node, when merging samplings.) \label{line:seed3d} \label{line:init3d}
		\FOR{$\vx\in X$}
		\STATE $G_{occ}\leftarrow G_{occ}\cup N^-(\vx)$ \hfill $\triangleright$ initialize occupied cells \label{line:hfhf}
		\ENDFOR

		\STATE	 \textbf{Sampling:} \label{line:sampling3d}
		\STATE The sampling process in 3D works exactly like in 2D with the two only difference being the following: 
		\begin{itemize}
			\item new candidates are generated on a spherical shell instead on an annulus
			\item a candidate $p\notin \delta D^3$ is rejected if $dist(p,\delta D^3)<\rho(p)/2$
		\end{itemize}
		\STATE \textbf{Resampling:}  (optional: algorithm terminates, if no resampling is required.)
		
		\REPEAT 
		\FOR{$T\in \mathcal{T}(X)$} \label{line:triangulate}
		\IF{$T$ is a sliver} \label{line:sliver}
		\STATE $X\leftarrow X\setminus\{\vx,\vy\}$, where $\vx,\vy\in T$ are 2 random nodes not contained in the boundary or the DFN \label{line:rem_nodes}
		
		 \hfill $\triangleright$ minimal distance of nodes to DFN and boundary assures, that this is possible.
		\ENDIF
		\ENDFOR

		\STATE \textbf{Rerun algorithm again from line \ref{line:sampling3d}}\label{line:rerun3d}
		\UNTIL{$\mathcal{T}(X)$ contains no more slivers.}\label{line:term3d}
		
	\end{algorithmic}

\end{algorithm}

\section{Results}\label{sec:results}

\subsection{Two-dimensional Examples}

\begin{figure}
	\centering
            \includegraphics[width=1\linewidth]{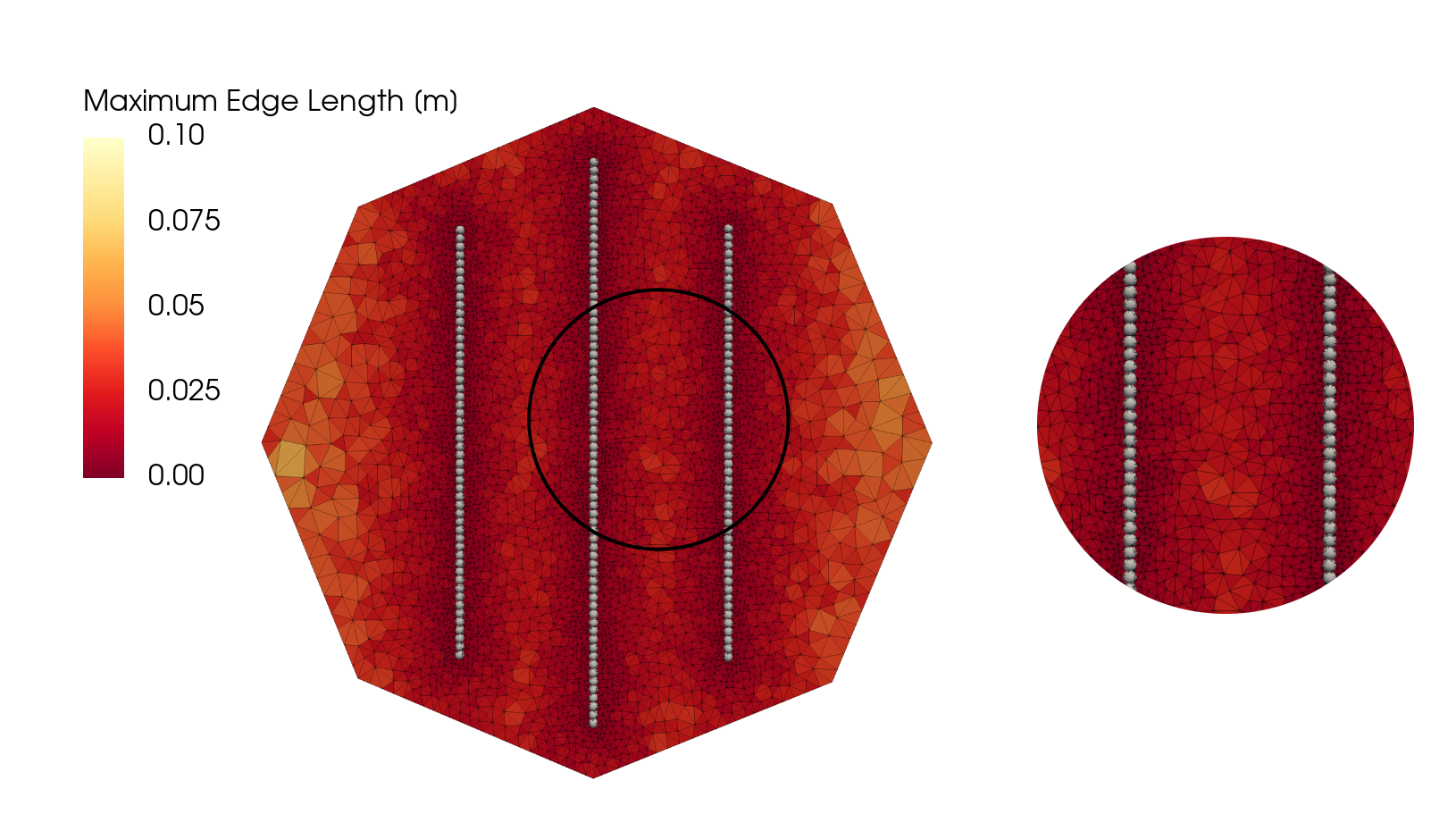}
	\caption{Triangulation of variable radii Poisson-disk sampling on fracture with three intersections. (H=0.01,R=40,A=0.1,F=1) Triangles colored according to their maximal edge length. The lines of intersection are shown as spheres. }
	\label{fig:varedgmaxfin}
\end{figure}

\begin{figure}
	\centering
	\includegraphics[height=0.4\linewidth]{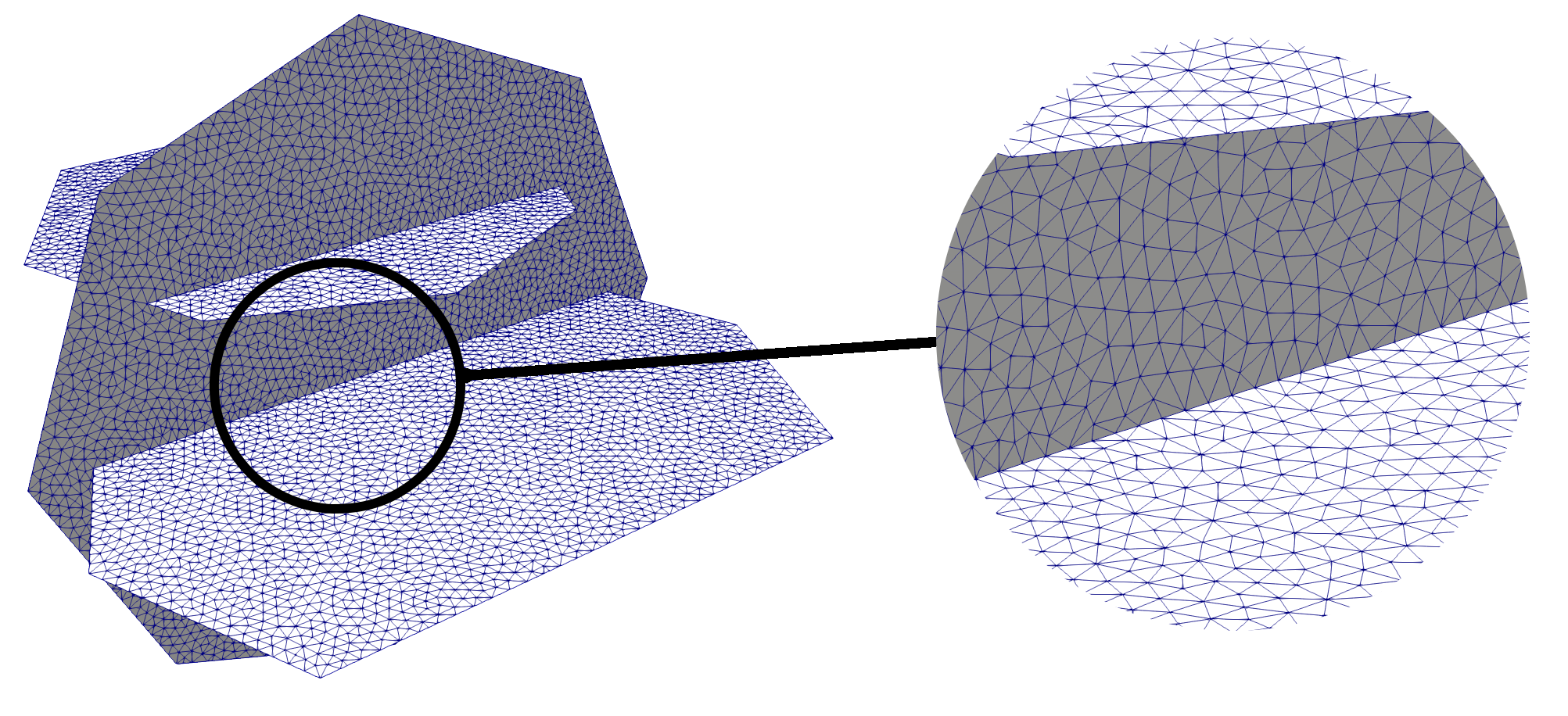}
	\caption{Triangulation of a regular Poisson-disk sampling reassembled into the original DFN.}
	\label{fig:2d-together}
\end{figure}

Figure \ref{fig:varedgmaxfin} shows the triangulation of a variable-radius sampling on a simple fracture with 3 intersections. Triangles are colored by their maximal edge length nicely showing how the triangle size increases as we move further away from the intersections.  
%
In Figure \ref{fig:2d-together}, we depict the triangulation of a constant-radius sampling on that same fracture, put back together into the original DFN it originated from. This process does not influence the overall triangulation quality unless the fractures themselves intersect in an angle smaller than the angles of triangles in the triangulation.  

\begin{figure}
	\centering
	\includegraphics[height=0.4\linewidth]{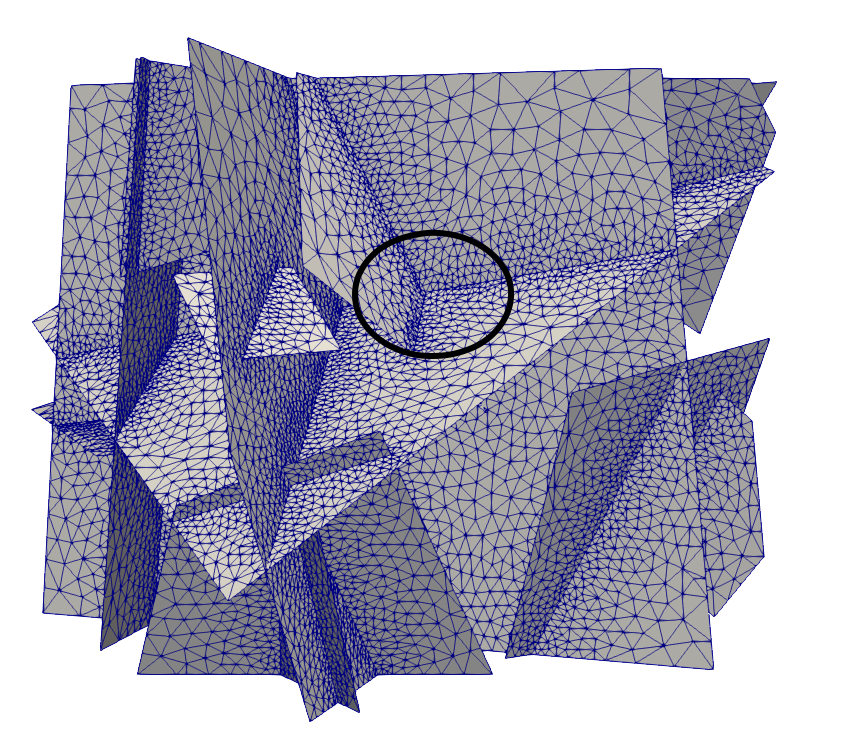} 
	\includegraphics[height=0.3\linewidth]{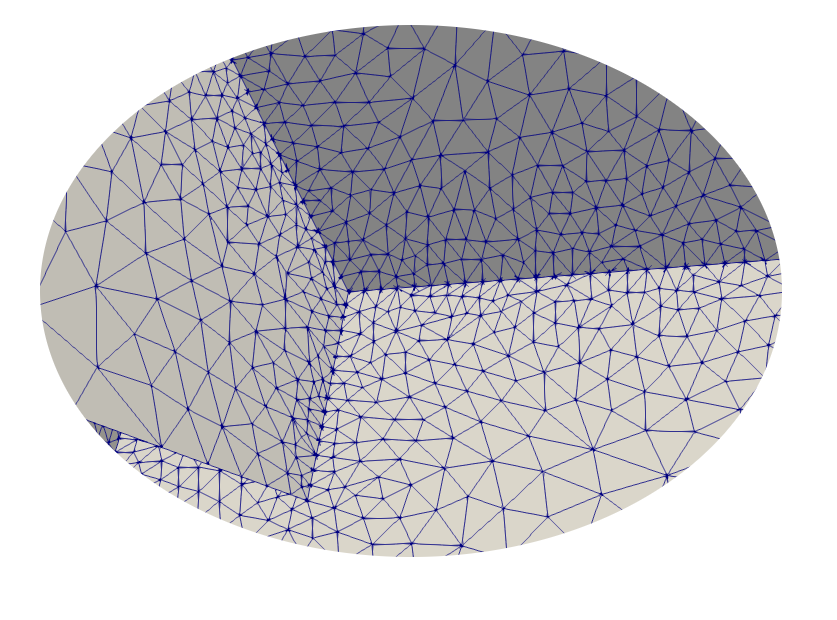}
	\caption{Triangulation of a variable radii Poisson-disk sampling reassembled into original DFN. 25 fractures, with up to eight intersections per fracture.($H=0.1$,$R=40$,$A=0.1$,$F=1$)}
	\label{fig:togetherzoomexp}
\end{figure}

\begin{figure}
	\centering
	\begin{tikzpicture}
	\node (A) at (0,0) {\input{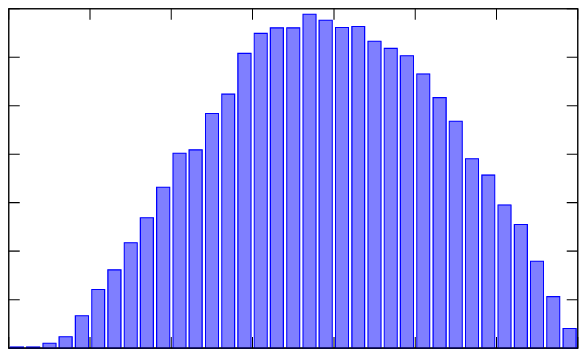}};
	\node (B) at (0,-5.5) {\input{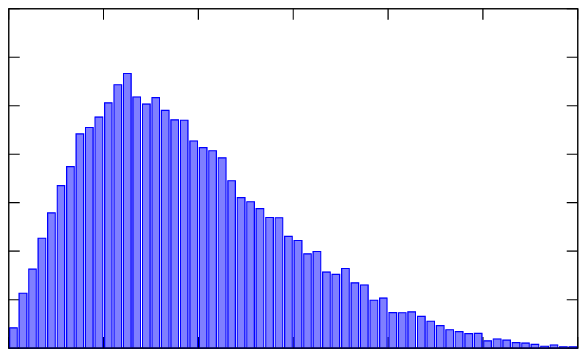}};
	\node (B) at (0,-11) {\input{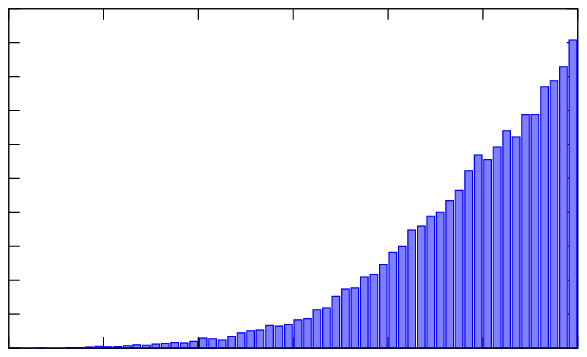}};
	\node (a) at (-4.,2.3) {(a)};
	\node (a) at (-4.,-3.2) {(b)};	
	\node (a) at (-4.,-8.7) {(c)};		
	\end{tikzpicture}

	\caption{Histograms of selected quality measures of the triangulation of variable radii Poisson-disk sampling on a fracture with three intersections. (H=0.01,R=40,A=0.1,F=1). (a): minimal angle ($\ge 25^\circ$), (b): max angle ($\le 120^\circ$), (c): aspect ratio ($\ge 0.47$)}
	\label{fig:variablemaxangle}
\end{figure}

We show an example from a slightly bigger DFN combining both variable radii Poisson-disk sampling and the reassembly into its original form in Figure \ref{fig:togetherzoomexp}. This particular example contains 25 fractures, with up to eight intersections per fracture, some of these intersecting each other. The parameters of the inhibition radius are set to $H=0.1$,$A=0.1$,$F=1$ and $R=40$.
\jeffrey{Add information about the network. Number of fractures, Length distribution, etc. }
The high quality of this particular triangulation is showcased in the histograms in Figure \ref{fig:variablemaxangle}. Depicted are the distribution of minimal angles (a), maximal angles and the aspect ratios of the triangulation. We see one triangle each with $25^\circ$ and $26^\circ$ respectively as minimal angles with all other minimal angles being greater than $27^\circ$. The theoretical minimum angle in a maximal Poisson-disk sampling with Lipschitz constant $A=0.1$ is $27.04^\circ$. The majority of minimal angles is significantly better still. In terms of the maximal angle, we can observe very few triangles with angles worse than $110^\circ$ and none worse than $120^\circ$. The greatest maximal angle theoretically possible in a maximal Poisson-disk sampling with this Lipschitz-constant would be $125.92^\circ$. The vast majority of aspect ratios is greater than $0.8$ with only a marginal number of triangles having an aspect ratio of less than $0.6$ and none below $0.47$.

\subsection{Run Time Analysis}
\begin{figure}
	\centering
	\input{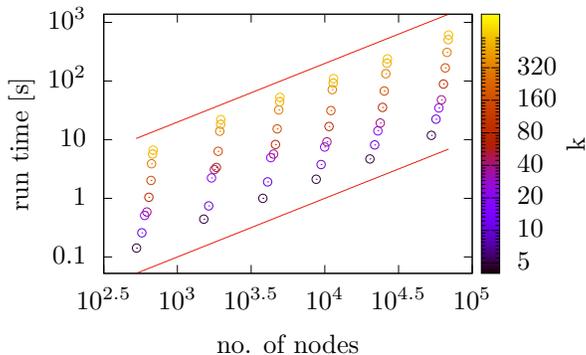}

	\caption{Log-log-plot of run time of Poisson-disk sampling algorithm in dependency of the total number of nodes sampled prior to the resampling process. Data points generated over the same DFN, different point densities generated by changing the minimal inhibition radius  $\frac{H}{2}$ between every pair of nodes. Data points are colored depending on the value of $k$. Other parameters are set to $A=0.1,R=40, F=1$. Comparison to lines of slope $1$ (red) indicates the run time increases approximately at a linear rate.	\label{fig:part-time1}}
\end{figure}

\begin{figure}
	\centering
	\input{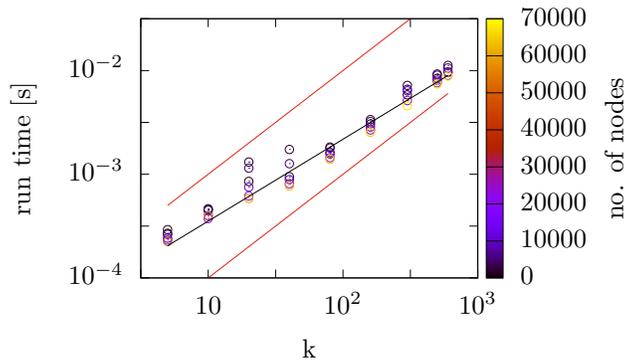}
	\caption{Log-log-plot of run time of Poisson-disk sampling algorithm in dependency of the number of concurrently sampled nodes $k$ prior to the resampling process. Data points generated over the same DFN, different point densities generated by changing the minimal inhibition radius  $\frac{H}{2}$ between every pair of nodes. Data points are colored depending on the total number of nodes sampled. Other parameters are set to $A=0.1,R=40, F=1$. Linear fit (black) with slope $0.7(9)\pm 0.00(7)$. Comparison to lines of slope $1$(red) indicate sublinear behavior.}
	\label{fig:ktime}
\end{figure}

We show an analysis of the run time and quality of the sampling on a DFN for varying sample sizes, variations of the parameter $k$, and different numbers of resampling attempts. All these data points were generated on the same DFN. Different node numbers were achieved by continuously changing the parameter $\frac{H}{2}$, the minimal allowed distance between nodes. All data points are from independent samplings.
The plot in Figure \ref{fig:part-time1} shows the run time prior to resampling process against the number of nodes sampled up to that point.
The color corresponds to the value of the parameter $k$, which controls the number of concurrent samples. 
We see an increase in run time with increasing $k$, as expected. The run times for samples with the same $k$ are positioned along straight lines of slope one, indicating a linear dependence of the total run time and the number of nodes sampled. The red lines in the plot have a slope of 1 to help visualize this. 
Figure \ref{fig:ktime} shows the relation between the parameter $k$ and the run time. Colors correspond to different numbers of nodes. As already established, the run time increases linearly with the number of nodes sampled.  The run time in terms of $k$ even exhibits a slightly sublinear behavior. The linear fit (black) of the data in this log-log-plot has a slope of $0.7(9)\pm 0.00(7)$. While this fitting error of $\approx 9\%$ is not insignificant it can also clearly be seen by comparing the data to the two lines of slope $1$ (red) in the plot, that the run time does not increase more than linearly with $k$.\\

\begin{figure}
	\centering
	\input{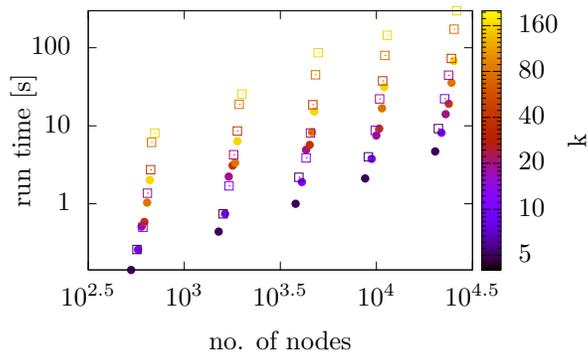}
	\caption{Comparison of run time for an implementation of \cite{Bridson2007FastPD,fastvar} (squares) and our variation ot the algorithm(circles). Depicted in a Log-log-plot are run time of Poisson-disk sampling algorithm in dependency of the total number of nodes sampled prior to the resampling process. Data points generated over the same DFN, different point densities generated by changing the minimal inhibition radius  $\frac{H}{2}$ between every pair of nodes. Data points are colored depending on the value of $k$. Other parameters are set to $A=0.1,R=40, F=1$. }
	\label{fig:nodestimecomp}
\end{figure}

Figure \ref{fig:nodestimecomp} depicts a comparison of runtime between our implementation of \cite{fastvar}  or \cite{Bridson2007FastPD} for variable-radii sampling and the same implementation with our adaptation to use the grid not only to find closeby nodes, but also directly reject candidates. Data points generated by our adapted algorithm are represented by a filled circle, whereas data points generated by the original algorithm are shown by empty squares. All data points are colored depending on $k$. We can see our algorithm out performs the original for every pair of data points. This advantage increases with growing $k$, which makes sense as there are more rejected candidates the greater $k$ is and our adapted version can handle rejection faster since it does not have to calculate the distance. For $k = 5$ the speed difference between the algorithms is slightly less than a factor of 2, whereas for $k = 160$ the advantage grows to about an order of magnitude. 
Comparisons to the sampling algorithm used in {\sc dfnWorks} \cite{hyman2015dfnWorks} prior to this implementation showed the current version is faster by yet another order of magnitude.

\subsection{Quality and resampling}
\begin{figure}
	\centering
	\input{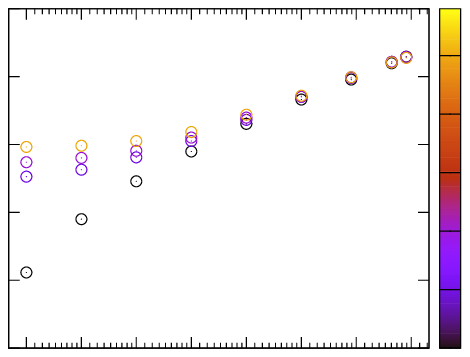}
	\caption{Total number of nodes sampled after resampling plotted in dependence of $k$ colored by number of resamplings. Data points generated over the same DFN with fixed minimal inhibition radius. Other parameters are set to $A=0.1,R=40, F=1$.}
	\label{fig:k_nodes_rep}
\end{figure}

The maximality of our samples correlates to a high degree to the choice of $k$, but also to the number of times the resampling algorithm is run. Depicted in Figure \ref{fig:k_nodes_rep} are the total number of nodes sampled after a different number of resamplings. First we can see that the density of nodes grows with the parameter $k$. This growth starts out fast for small $k$ and while not entirely ceasing to increase, slows down notably for higher $k$. (Note log-scale on $\vx$-
axis.) On the lower end of the $k$ scale, resampling increases the node density significantly, whereas there is barely any difference for higher $k>100$. 
The first resampling is particularly effective, whereas the difference between each resampling decreases afterwards. Given that resampling does not take more time then the original sampling process, this turns into an interesting trade-off between higher $k$ and more repetitions of the resampling that overall can yield higher performance. A run at $k=5$ with few repetition for example, results in a density comparable to a run with more than $10$ times higher $k$ without resampling, while being significantly faster overall. Similar conclusions can be reached when looking at the quality of resulting triangulations rather than just the density of the Poisson-disk sampling.

\begin{figure}
	\centering
	\input{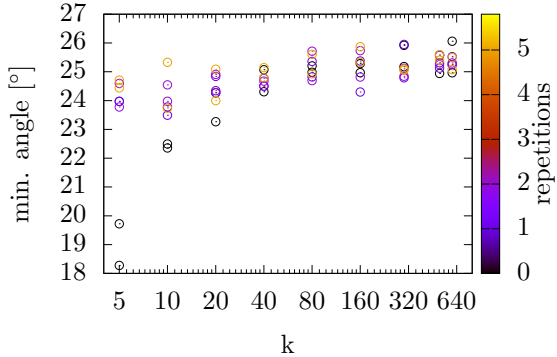}
	\caption{Smallest minimal angle of the triangulation of Poisson-disk samplings in dependence of $k$ colored by the number of resamplings. Data points generated over the same DFN with fixed minimal inhibition radius.   Other parameters are set to $A=0.1,R=40, F=1$.}
	\label{fig:k_ang}
\end{figure}

Figure \ref{fig:k_ang} shows the smallest minimal angle in a triangulation of our sampling for variable $k$ and different numbers of resampling attempts. We can see for $k\gtrapprox 80$ this angle appears to be at around $25^\circ$ independently of the number of repetitions. 
The theoretical bound for a maximal Poisson-disk sampling (with $r(\vx,\vy)=R(\vx,\vy )$) for the settings used to generate these data points would be $27.04^\circ$. Solving the the angle bounds from lemma \ref{lemma:2} for $\varepsilon$ shows us that in this sampling $R(\vx,\vy)\lesssim (1+0.1)r(\vx,\vy)$.
Given the statistical nature of the algorithm and the fact that identical inhibition and coverage radii are not quite guaranteed these results can be considered very good.
While the quality of triangulations for smaller $k$ without resampling is significantly lower, it is noteworthy that just a single repetition fixes this issue and yields triangulations with qualities on par with those for even significantly higher $k$. 
This allows the algorithm to run at single or low double digit $k$, perform a single resampling and generate a triangulation just as good as higher $k$ would have produced in multitudes of the time.

\subsection{Three-Dimensional Example}

\begin{figure*}
	\centering
\includegraphics[height=0.3\textwidth]{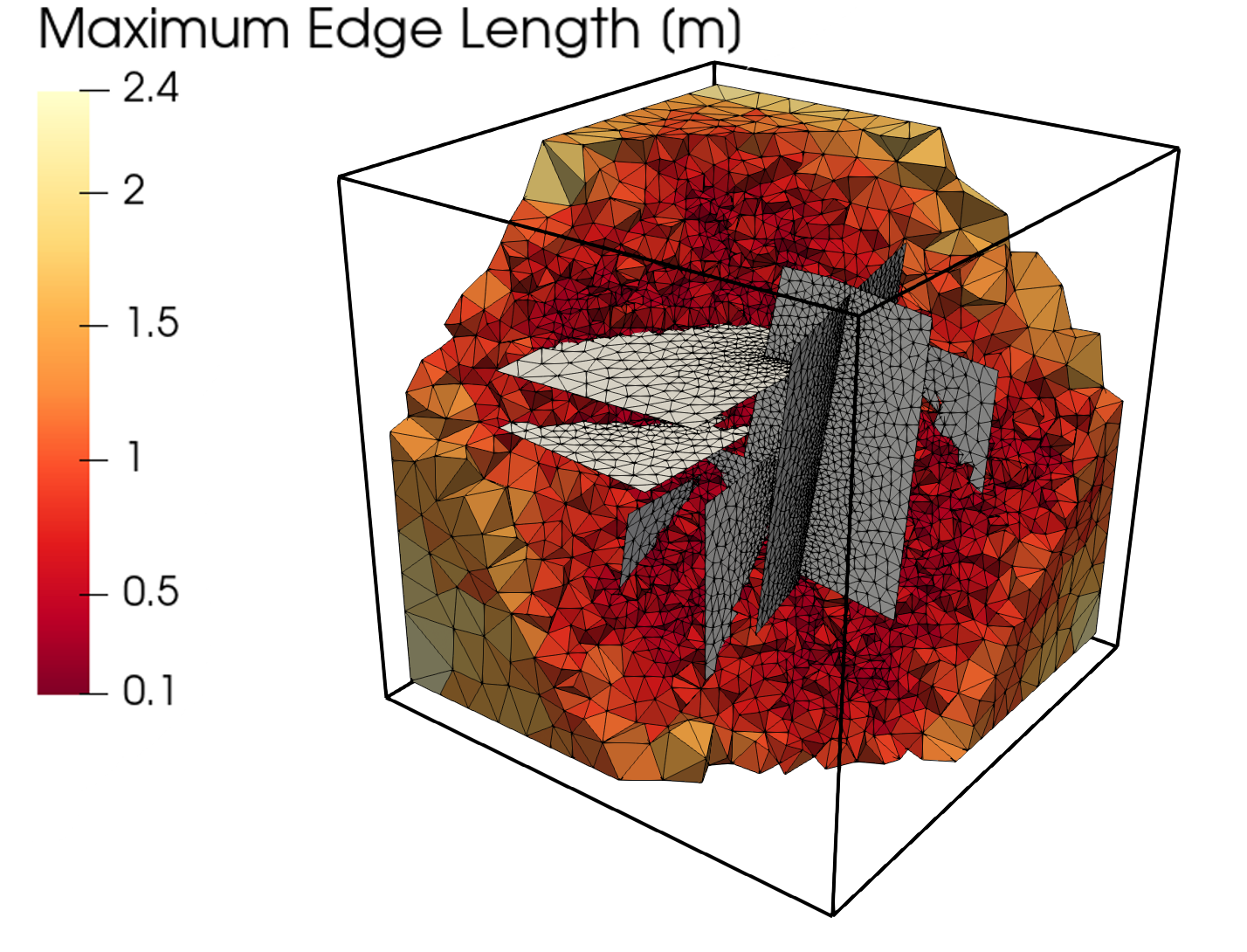}	\includegraphics[height=0.3\textwidth]{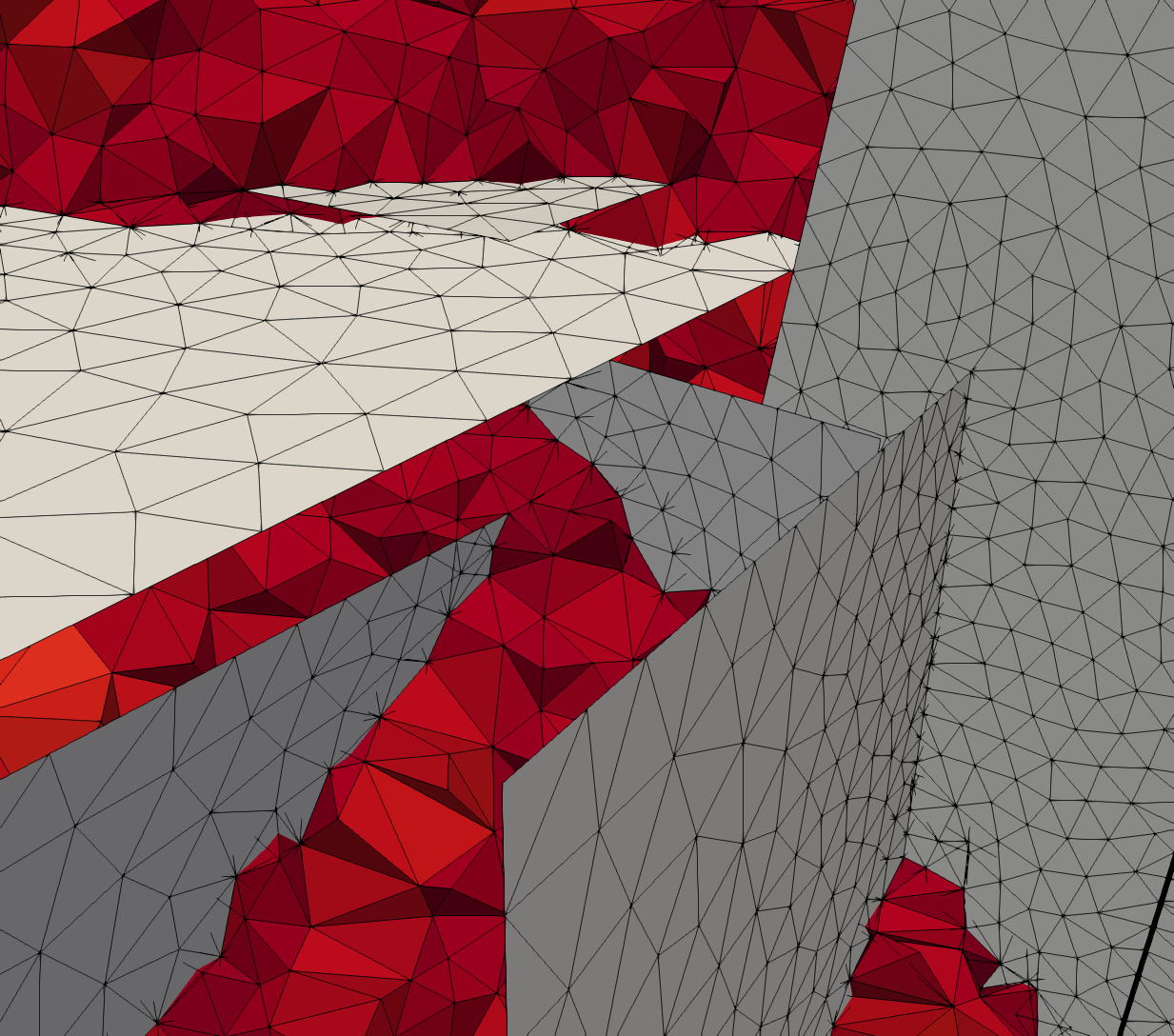}
	\caption{(Left) Triangulation of variable radii Poisson-disk sampling of DFN and its surrounding region. (Right) Close up of the conforming mesh. (H=0.25,R=100,A=0.125,F=1). Tetrahedra colored according to their maximal edge length.}
	\label{fig:3d-mixzoomfin}
\end{figure*}

While the majority of our work was aimed at optimizing the 2D sampling on a DFN, we will conclude with an example where these 2D samplings are combined with a 3D sampling of the surrounding matrix to showcase that it can be used to produce high quality triangulations in this case as well. Triangulated output of the 3D algorithm can be seen in Figure \ref{fig:3d-mixzoomfin}. The tetradedra are colored according to their maximal edge length to show how the point density is adapted with the distance to the closest fracture.

\begin{figure}
	\centering
	\begin{tikzpicture}
	\node (A) at (0,0) {	\input{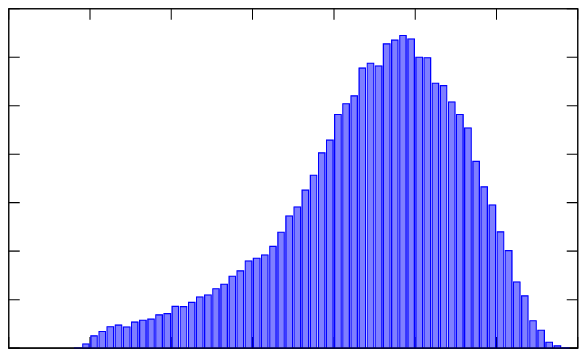}};
	\node (B) at (0,-5.5) {	\input{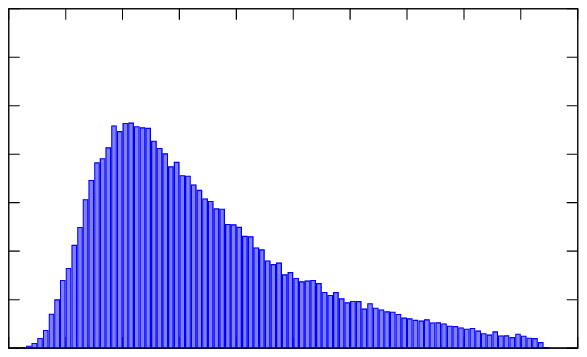}};
	\node (B) at (0,-11) {\input{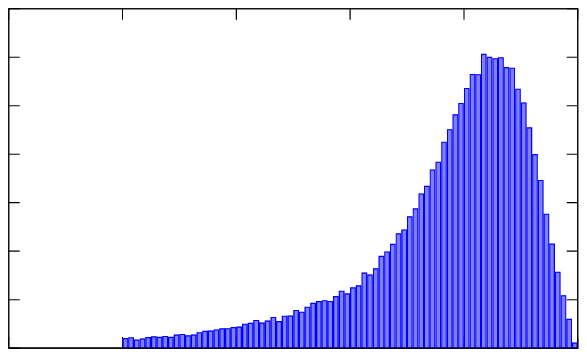}};
	\node (a) at (-4.,2.3) {(a)};
	\node (a) at (-4.,-3.2) {(b)};	
	\node (a) at (-4.,-8.7) {(c)};		
	\end{tikzpicture}
	
	\caption{Histograms of selected quality measures of the triangulation of variable radii Poisson-disk sampling on DFN and its surrounding matrix. (H=0.01,R=40,A=0.1,F=1). (a): minimal angle ($\ge 8^\circ$), (b): max angle ($\le 165^\circ$), (c): aspect ratio ($\ge 0.2$)}
	\label{fig:3dhist}
\end{figure}

Finally, the histograms in Figure \ref{fig:3dhist} show the distribution of quality measures of the tetrahedra in the triangulation depicted in Figure \ref{fig:3d-mixzoomfin}. 
For this run, tetrahedra with either a dihedral angle of less than $8^\circ$ or an aspect ratio of less than $0.2$ were discarded before the sampling algorithm was restarted.  
The first histogram depicts the distribution of the minimal dihedral angle of each tetrahedron. 
As expected no dihedral angle below $8^\circ$ remains, while the vast majority exceeds values of $30^\circ$. 
Histogram (b) shows that despite not optimizing with respect to the maximal dihedral angle none of these angles exceed $165^\circ$. 
Histogram (c) shows a sharp cut-off at $0.2$ in the distribution of aspect ratios indicating that the aspect ratio is likely to have been  the driving factor for a majority of the resamplings.  
The example shown ran through a sliver-removal and resampling process $17$ times to obtain its triangulation quality. In each of these steps a total of $200$ or less out of approximately $50000$ nodes were removed before the resampling.

\section{Conclusions}\label{sec:conclusions}
We presented algorithms that successfully generate variable-radii Poisson-disk samples on polygonal regions or networks of polygons and the surrounding space they are embedded in. We increased the performance of existing algorithms and introduced additional measures to guarantee certain levels of maximality. 
It is worth noting that maximality is reached for a coverage radius just slightly larger than the inhibition radius.
Triangulations of these samplings show a quality almost matching theoretical quality bounds for maximal Poisson-disk samplings, in which coverage and inhibition radii coincide.
Our key contributions are summarized as:
\begin{enumerate}
	\item our algorithm is significantly faster than the previous conforming variable mesh strategies
	\item for the fracture networks, we achieved mesh quality only marginally worse than what is theoretically possible,
	\item for the volume meshing, slivers can be removed entirely from the domain within certain bounds
\end{enumerate}  


It is worthwhile mentioning that the described algorithms are  not only fast, but also simple to run in a parallel fashion, further improving the overall runtime. 
Given a DFN, the 2D-sampling can parallelized by working on each fracture on a different processor. 
Based on the grid structure used to accept and reject candidates, both 2D and 3D can also be further parallelized by dividing their domain into several pieces, that can be sampled individually on different processors, while needing to communicate only cell  information on the boundaries of the split domains. 
Once these point distributions are produced, however, the all must reside on a single processor to connect them into a Delaunay mesh.

\section{Acknowledgments}
J.K. gratefully acknowledges support from the 2020 National Science Foundation
Mathematical Sciences Graduate Internship to conduct this research at Los Alamos National Laboratory.
J.D.H. and M.R.S. gratefully acknowledges support from the LANL LDRD program office Grant Number \#20180621ECR, the Department of Energy Basic Energy Sciences program (LANLE3W1), and the Spent Fuel and Waste Science and Technology Campaign, Office of Nuclear Energy, of the U.S. Department of Energy.
M.R.S. would also like to thank support from the Center for Nonlinear Studies.
J.M.R.  received support from DOE,  Contract No. DE-AC05-00OR22725. 
Los Alamos National Laboratory is operated by Triad National Security, LLC, for the National Nuclear Security Administration of U.S. Department of Energy (Contract No. 89233218CNA000001).
This work was prepared as an account of work sponsored by an agency of the United States Government.
The views and opinions of authors expressed herein do not necessarily state or reflect those of the United States Government or any agency thereof, its contractors or subcontractors.
LAUR \# LA-UR-21-24804.





	\bibliography{references}{}
\bibliographystyle{plain}
\end{document}